\documentclass[AMS,STIX1COL]{WileyNJD-v2}
\articletype{Article Type}

\newboolean{boldeqnitalic}
\setboolean{boldeqnitalic}{true}
\usepackage{defdb1}
\usepackage{soul}
\received{5 July 2023}
\revised{X XXXX 2023}
\accepted{X XXXX 2023}

\usepackage{caption}
\usepackage{subcaption}
\usepackage{cleveref}
\crefname{equation}{}{}  

\def\Cpp{{C\nolinebreak[4]\texttt{++}}}

\begin{document}

\title{
A computational framework for pharmaco-mechanical interactions in arterial walls using parallel monolithic domain decomposition methods}

\author[1]{D. Balzani}
\author[2]{A. Heinlein}
\author[3,4]{A. Klawonn}
\author[3,4]{J. Knepper}
\author[1]{S. Nurani Ramesh}
\author[5,6]{O. Rheinbach}
\author[3]{L. Sa\ss mannshausen}
\author[1]{K. Uhlmann}

\address[1]{\orgdiv{Civil and Environmental Engineering}, \orgname{Ruhr University Bochum}, \orgaddress{Universitätsstraße 150, 44801 Bochum, Germany}}
\address[2]{\orgname{Delft Institute of Applied Mathematics}, \orgaddress{Delft University of Technology, Mekelweg 4, 2628 CD Delft, The Netherlands}}
\address[3]{\orgname{University of Cologne}, \orgaddress{Department of Mathematics and Computer Science, Cologne, Germany}}
\address[4]{\orgname{University of Cologne}, \orgaddress{Center for Data and Simulation Science, Cologne, Germany}}
\address[5]{\orgname{Technische Universität Bergakademie Freiberg}, \orgaddress{Fakultät für Mathematik und Informatik,
Freiberg, Germany}}
\address[6]{\orgname{Technische Universität Bergakademie Freiberg}, \orgaddress{Zentrum für effizente Hochtemperaturstoffwandlung (ZeHS), Freiberg, Germany}}

\corres{Axel Klawonn, University of Cologne, Department of Mathematics and Computer Science, Weyertal 86-90, 50931 K\"oln, Germany. \email{axel.klawonn@uni-koeln.de}}

\abstract[Abstract]{
A computational framework is presented to numerically simulate the effects of antihypertensive drugs, in particular calcium channel blockers, on the mechanical response of arterial walls. 
A stretch-dependent smooth muscle model by Uhlmann and Balzani is modified to describe the interaction of pharmacological drugs and the inhibition of smooth muscle activation. 
The coupled deformation-diffusion problem is then solved using the finite element software FEDDLib and overlapping Schwarz preconditioners from the Trilinos package FROSch. 
These preconditioners include highly scalable parallel GDSW (generalized Dryja–Smith–Widlund) and RGDSW (reduced GDSW) preconditioners. 
Simulation results show the expected increase in the lumen diameter of an idealized artery due to the drug-induced reduction of smooth muscle contraction, as well as a decrease in the rate of arterial contraction in the presence of calcium channel blockers. Strong and weak parallel scalability of the resulting computational implementation are also analyzed.
}


\maketitle

\section{Introduction}
\label{sec:intro}

Approximately one in four adults worldwide has hypertension \cite{forouzanfar2017global}. In contemporary societies, the mean blood pressure levels increase steadily with age, in stark contrast to pre-industrial societies, where the levels changed little with age \cite{oparil2018hypertension}. Furthermore, there has been a rapid increase in the prevalence of hypertension in the last few decades \cite{mills2016global}. Hypertension is a major risk factor for cardiovascular disease and is associated with atherosclerosis, stroke, heart failure, and kidney disease. Atherosclerosis is a condition characterized by a reduction in the arterial lumen diameter and obstructed blood flow, due to the formation of fibrofatty lesions called atherosclerotic plaques. These plaques contain oxidized lipids, calcium deposits, inflammatory cells, and smooth muscle cells.

Hypertension and atherosclerosis are the two major causes of cardiovascular disease, which is the most common form of mortality in the world. Both these conditions are treated by antihypertensive drugs such as ACE inhibitors, angiotensin II receptor blockers, dihydropyridine calcium channel blockers and thiazide diuretics. Some of these drugs directly affect the arterial wall and can lead to changes in wall structure and function. For example, calcium channel blockers work by reducing the contractility of the arterial wall, thereby reducing blood pressure. However, it has been observed that when multiple drugs are used to treat hypertension, the risk of ischemic events increases with each drug class. In cases where three or more classes of drugs were necessary to successfully control hypertension, the risk of stroke was observed to be 2.5~times higher than in healthy normotensive people~\cite{howard2015blood}. Additionally, the exact mechanism of the drug-induced mechanical response of the arterial wall in atherosclerotic arteries is not well understood. In light of the increasing prevalence of hypertension and the resulting uptake of antihypertensive drugs, there is a pressing need for a better understanding of pharmaco-mechanical interactions in arteries. To this end, numerical simulation of healthy and atherosclerotic arteries can be a valuable tool. Simulating the mechanical response entails an accurate description of the arterial wall material behavior, the drug transport process and the drug-artery interaction.

Arteries consist of three major components: endothelial cells, smooth muscle cells (SMCs) and the extracellular matrix (ECM). The thickness of the endothelium is negligible in comparison to the thickness of the artery and so is its load-carrying capacity; therefore it is not considered in our model. The ECM is composed of elastin and collagen fibers and supports the mechanical load. Vascular SMCs are arranged concentrically and their coordinated contraction and relaxation causes changes in luminal diameter. Contraction in SMCs is regulated by cytosolic free calcium concentration. A change in membrane potential due to factors like cell stretch, neurotransmitters, and hormones causes the opening of voltage-gated calcium channels, thereby enabling the movement of calcium ions into the cytosol \cite{touyz2019textbook}. The calcium forms a complex with calmodulin and activates the enzyme myosin light chain kinase (MLCK), which phosphorylates the regulatory light chain (RLC) of myosin, enabling contraction. The enzyme myosin light chain phosphatase (MLCP) is responsible for the dephosphorylation of myosin RLC. Enzymes such as Rho kinase inhibit the activity of MLCP and lead to a calcium-independent contraction mechanism \cite{webb2003smooth}. Since antihypertensive drugs primarily work by interacting with SMC activation, it is important to have an SMC model that allows for a meaningful description of drug-tissue interaction. There are several models that describe the active response of arteries \cite{murtada2010calcium,murtada2012experiments,yang2003myogenica,yang2003myogenicb,bol2012three,staalhand2011mechanochemical,uhlmann2023chemo}. The well-accepted cross-bridge phosphorylation model by Hai and Murphy~\cite{hai1988cross} describes the influence of MLCK and MLCP on contraction. Yang \etal~\cite{yang2003myogenica,yang2003myogenicb} proposed an electro-chemo-mechanical model where the change in membrane potential due to the various ion channels is considered. In B\"ol \etal~\cite{bol2012three}, the calcium concentration was defined as a function of time, and arteries with calcium waves were simulated with a one-sided coupling. Uhlmann and Balzani \cite{uhlmann2023chemo} extended the model proposed by Murtada \etal~\cite{murtada2012experiments} to include the influence of stretch on the action of MLCP and MLCK. Here, we adopt the active response model by Uhlmann and Balzani \cite{uhlmann2023chemo} and modify it to model the impact of dihydropyridine calcium channel blockers (CCBs). In a previous work \cite{nurani_ramesh_first_2023}, we already investigated the effect of CCBs on the activity of MLCK, considering constant MLCP activity. CCBs work by blocking the voltage-gated calcium channels and thereby reducing the calcium influx into the cytosol. The arterial wall mass transport process, depending on the type of drug, is either advection or diffusion dominated and can therefore be modeled by a linear scalar advection-diffusion equation \cite{ai2006coupling,hossain2012mathematical,hwang2001physiological}. However, due to drug absorption and binding, and to account for complex patient-specific arteries, the three-dimensional reaction-advection-diffusion equation is necessary to comprehensively describe the drug transport phenomenon in arteries. In the case of hydrophobic drugs, diffusion-dominated transport is observed, whereas, in hydrophilic drugs, advection-dominated transport is observed \cite{hwang2003impact,creel2000arterial,o2010factors}. Since we consider only hydrophobic CCBs, we may simplify the transport process and model it using the reaction-diffusion equation. 

To simulate patient-specific arteries, a sufficient spatial resolution is required, leading to large number of degrees of freedom and ill-conditioned matrices. 
Moreover, for the interaction of the fluid with the arterial wall, strong coupling schemes are most suitable, and monolithic fluid-structure interaction (FSI) coupling schemes are most competitive in this regime; this behavior can be mathematically explained by the so-called \emph{added-mass effect}~\cite{causin_added-mass_2005}. This effect also explains instabilities typically resulting from the use of weakly coupled schemes, which are very efficient in other applications like aeroelasticity.
	
In this work, the effect of the antihypertensive drug in the structure is modeled and simulated; i.e., we are concerned with coupled structure-chemistry interaction problem (SCI). A fully monolithic approach is taken; i.e., the multiphysics problem is assembled into a single system.
To solve this problem, we then require robust parallel scalable preconditioners. 

The remainder of the paper is organized as follows. In Section 2, the model for the arterial wall is introduced and in Section 3, the pharmaco-mechanical interaction is modeled. Next, in Section 4, the monolithic, two-level overlapping Schwarz preconditioners for the coupled system are presented. In Section 5, the software ecosystem defining our computational framework is described. Then, numerical results obtained using this software framework are presented. In Section 7, a conclusion is given.

\section{Modeling the arterial wall} \label{sec:modeling}
In the present work, we restrict ourselves to resistance arteries, which are characterized by a high SMC content. Since we are interested in the vasoregulatory action of arteries, which is primarily a stretch-dependent, chemo-mechanically coupled problem, we neglect SMC activation due to other factors like the sympathetic nervous system and the endocrine system. We further restrict our model to the tunica media layer of the artery, which consists of a passive extracellular matrix (ECM) and active smooth muscle cells. The passive behavior of the ECM is modeled using a hyperelastic material law described in Balzani \etal~\cite{balzani2006polyconvex}. Two distinct crosswise helically arranged fiber directions are considered, both lying in the longitudinal-circumferential plane and symmetric about the circumferential axis. The active material is modeled using the phenomenological model by Uhlmann and Balzani~\cite{uhlmann2023chemo}, where the models of Hai and Murphy~\cite{hai1988cross} and Murtada \etal~\cite{murtada2012experiments} are extended to include the effects of stretch on the calcium-dependent and -independent contraction mechanisms.

\subsection{General continuum-mechanical model}
Let $\bX$ and $\bx$ denote a material point in the initial configuration $\mathcal{B}$ and the deformation in the current configuration $\mathcal{S}$, respectively. The motion of the point over time is defined by the map $\bx = \varphi(\bX, t)$. The associated deformation gradient $\bF$ and the right Cauchy--Green tensor $\bC$ are defined as
\begin{equation}
    \bF = \frac{\partial \bx}{\partial \bX}, \quad \bC = \bF^T \bF.
\end{equation}
We adopt an invariant-based framework for the constitutive model with two discrete fiber directions $\ba^{(f)}$, where the strain energy density $\Psi$ is additively decomposed as
\begin{equation}
    \Psi = \Psi_{\mathrm{p,isot}} + \sum_{f=1}^{2} \Psi_{\mathrm{p,ti}}^{(f)} + \sum_{f=1}^{2} \Psi_{\mathrm{a}}^{(f)},
\end{equation}
since it is assumed that there is a weak interaction between the fiber families. $\Psi_{\mathrm{p,isot}}$ represents the influence of the isotropic elastin matrix, $\Psi_{\mathrm{p,ti}}^{(f)}$ - the transversely isotropic part models the influence of the passive collagen fibers, whereas $\Psi_{\mathrm{a}}^{(f)}$ represents the response of the active smooth muscle cells. The principal and mixed invariants of the right Cauchy--Green tensor $\bC$ are defined as
\begin{equation}
        I_1 = \tr(\bC), \quad I_2 = \frac{1}{2} \left(\tr(\bC)^2 - \tr(\bC^2)\right), \quad I_3 = \det(\bC), \quad I_4^{(f)} = \bC : \bM^{(f)}, \quad I_5^{(f)} = \bC^2 : \bM^{(f)},
\end{equation}
where $\bM^{(f)} = \ba^{(f)} \otimes \ba^{(f)}$ are the structural tensors. Using the formulation by Balzani \etal~\cite{balzani2006polyconvex} we have the passive components
\begin{equation}
    \begin{aligned}
        \Psi_{\mathrm{p,isot}} = \alpha_1\left(I_1 I_3^{-1/3} - 3\right)
        &+ \alpha_2\left(I_3^{\alpha_3} + I_3^{-\alpha_3} - 2\right), \\
        \Psi_{\mathrm{p,ti}}^{(f)} &= \alpha_4 \langle K_3^{(f)} - 2 \rangle^{\alpha_5},
    \end{aligned}
\end{equation}
where $\alpha_1$, $\alpha_2$, $\alpha_3$, $\alpha_4$ and $\alpha_5$ are material parameters and $K_3^{(f)} = I_1 I_4^{(f)} - I_5^{(f)}$.

\subsection{Stretch-dependent chemical kinetics model for smooth muscle cells} 
\label{sec:ChemicalKinetics}
The SMC cytoskeleton is composed of a three-dimensional network of actin and myosin filaments connected at cytoplasmic dense bodies. SMC contraction is a result of relative sliding between the myosin and actin filaments. The myosin filaments contain distinct head regions that interact with actin to form attached cross-bridges. Contraction primarily occurs due to increased intracellular free calcium concentration $\left(\left[\text{Ca}^{2+}\right]_\textrm{i}\right)$. Calmodulin, a calcium-binding protein, binds with calcium ions to form a complex that activates the enzyme myosin light chain kinase (MLCK), which results in the phosphorylation of the 20-kDa light chain of myosin. This stimulates myosin-ATPase activity, enabling cross-bridge cycling and force generation, resulting in SMC contraction. The myosin regulatory light chain is dephosphorylated by the activity of the enzyme myosin light chain phosphatase (MLCP). Hai and Murphy \cite{hai1988cross} proposed a model where cross-bridges exist in four distinct states: free dephosphorylated (A), free phosphorylated (B), attached phosphorylated (C) and attached dephosphorylated (D). Their change is described by seven rate constants: $k_1$ and $k_6$ represent the rate of phosphorylation by MLCK, $k_2$ and $k_5$ represent the rate of dephosphorylation by MLCP, $k_3$ is the rate of cross-bridge attachment, $k_4$ and $k_7$ are the rates of cross-bridge detachment. The kinetic model is thus governed by the following system of ordinary differential equations
\begin{equation}
    \begin{bmatrix}\dot{n}_\text{A}\\ \dot{n}_\text{B}\\ \dot{n}_\text{C}\\ \dot{n}_\text{D}\end{bmatrix}=\begin{bmatrix}-k_1&k_2&0&k_7\\ k_1&-k_2-k_3&k_4&0\\ 0&k_3&-k_4-k_5&k_6\\ 0&0&k_5&-k_6-k_7\end{bmatrix}\begin{bmatrix}n_\text{A}\\ n_\text{B}\\ n_\text{C}\\ n_\text{D}\end{bmatrix},
\end{equation}
where $n_\text{A}$, $n_\text{B}$, $n_\text{C}$ and $n_\text{D}$ are the fractions of cross-bridges in the four aforementioned states. The rate of cross-bridge attachment ($k_3$) and detachment ($k_4$ and $k_7$) are assumed to be constant, whereas the rates of phosphorylation $k_1$ and $k_6$ are given, as suggested first by Murtada \etal~\cite{murtada2012experiments}, by 
\begin{equation}
    k_{1 / 6}=\eta \frac{\left[\mathrm{Ca}^{2+}\right]_{\text{i}}^2}{\left[\mathrm{Ca}^{2+}\right]_{\text{i}}^2+\left(\mathrm{Ca}_{50}\right)^2},
\end{equation}
where $\eta$ is a material parameter and $\mathrm{Ca}_{50}$ is the half activation constant. Further, as proposed by Uhlmann and Balzani~\cite{uhlmann2023chemo}, the cytosolic calcium concentration is modeled as a stretch-dependent function
\begin{equation}
    \label{eq:Ca}
    \left[\mathrm{Ca}^{2+}\right]_{\text{i}} = \gamma_1 \langle \lambda - \bar{\lambda}_c \rangle^2,
\end{equation}
with $\lambda=(I_4)^{1/2}$ being the stretch along the SMC longitudinal direction, and $\gamma_1$ a material parameter. The authors modeled the reduction in intracellular calcium concentration due to various factors like calcium sequestration and calcium outflow by an evolution equation for $\bar{\lambda}_c$
\begin{equation}
    \begin{aligned}
    \dot{\bar{\lambda}}_{\textrm{c}} =
    \dot{\bar{\lambda}}_{\textrm{c,min}} +
    \dfrac{\dot{\bar{\lambda}}_{\textrm{c,max}} - \dot{\bar{\lambda}}_{\textrm{c,min}}}
    {1 + e^{\gamma_2\left([\mathrm{Ca}^{2+}]_{\textrm{tar}} - [\mathrm{Ca}^{2+}]_{\textrm{i}} - \tau_{\textrm{c}} \right)}},
    \end{aligned}
\end{equation}
where $\dot{\bar{\lambda}}_{\textrm{c,min}}$, $\dot{\bar{\lambda}}_{\textrm{c,max}}$, $\gamma_2$, $\tau_{\textrm{c}}$ are parameters and $[\mathrm{Ca}^{2+}]_{\textrm{tar}}$, the target calcium concentration, is defined by a stretch-dependent function 
\begin{equation}
    \label{eq:CaTar}
    [\mathrm{Ca}^{2+}]_{\textrm{tar}} = \gamma_3 \dfrac{\lambda^2} {\lambda^2 + \lambda_{\textrm{50,c}}^2},
\end{equation}
with $\gamma_3$ being a material parameter and $\lambda_{\textrm{50,c}}$ the half activation stretch. Calcium sensitization, which is a calcium-independent contraction mechanism \cite{somlyo2003ca2+}, occurs as a result of the action of the enzymes RhoA and Rho-Kinase (ROK). Uhlmann and Balzani \cite{uhlmann2022simulation2} modeled the resulting reduction in the rate of dephosphorylation by evolution equations for $k_{2/5}$
\begin{equation}
    \label{eq:k25}
    \begin{aligned}
        \dot{k}_{2 / 5} =
        \dot{k}_{2 / 5, \min }\left(1-e^{-\zeta_1 \left(k_{2 / 5} - k_{2 / 5, \min}\right)}\right)+
        \dfrac{\dot{k}_{2 / 5, \max }-\dot{k}_{2 / 5, \min }}
        {1+e^{\gamma_4 \left( \lambda-\bar{\lambda}_{\mathrm{p}} -\tau_{\mathrm{p}}\right)}}, \\
        \dot{\bar{\lambda}}_{\mathrm{p}} = 
        \dot{\bar{\lambda}}_{\mathrm{p}, \min } +
        \dfrac{\dot{\bar{\lambda}}_{\mathrm{p}, \max }-\dot{\bar{\lambda}}_{\mathrm{p}, \min }}
        {1+e^{\gamma_5\left(k_{2 / 5,\mathrm{tar}} - k_{2/5}-\tau_{\mathrm{k}}\right)}} -
        \dot{\bar{\lambda}}_{\mathrm{p}, \max } e^{-\zeta_2\left(\lambda - \bar{\lambda}_{\mathrm{p}}-\Delta \bar{\lambda}_{\mathrm{p}, \min }\right)},
    \end{aligned}
\end{equation}
where $\dot{k}_{2 / 5,\min}$ and $\dot{k}_{2 / 5,\max}$ are the minimum and maximum rate of change of $k_{2/5}$ respectively, and $\gamma_4$ and $\tau_{\textrm{p}}$ are parameters. Here, we additionally use $k_{2/5, \min}$ along with the penalty parameter $\zeta_1$ to ensure a minimum rate of dephosphorylation. The time adaption of $k_{2/5}$ after a stretch change is given by the second part of Equation \eqref{eq:k25}, wherein $\dot{\bar{\lambda}}_{\mathrm{p,max}}$, $\dot{\bar{\lambda}}_{\mathrm{p,min}}$, $\gamma_5$ are parameters and $\zeta_2$ is a penalty parameter which ensures a slow relaxation of the SMC by limiting $\lambda - \lambda_{\mathrm{p}}$ to a minimum value of $\Delta \lambda_{\mathrm{p,min}}$. Further, the target stretch-dependent MLCP activity is given as
\begin{equation}
    k_{2 / 5,\mathrm{tar}} = \gamma_6 \left(1 - \dfrac{\lambda}{\lambda + \lambda_{\mathrm{50,p}}}\right),
\end{equation}
where $\gamma_6$ is a material parameter and $\lambda_{\mathrm{50,p}}$ is the half activation stretch.

\subsection{Effect of pharmacological agents}
Here, we are interested in modeling the effects of calcium channel blockers (CCBs) on intracellular calcium concentration. CCBs typically reduce the flow of calcium ions through L-type calcium channels located on the smooth muscle cell membrane. Clinically, dihydropyridine CCBs are effective vasodilators and are widely used as antihypertensives. Their effect can be modeled by extending the chemical kinetics model given in \Cref{sec:ChemicalKinetics} by making the target calcium concentration $[\mathrm{Ca}^{2+}]_{\textrm{tar}}$ and the intracellular calcium concentration $[\mathrm{Ca}^{2+}]_{\textrm{i}}$ dependent on the concentration of the CCB. This is accomplished by making the parameters $\gamma_1$ and $\gamma_3$ in Equations \eqref{eq:Ca} and \eqref{eq:CaTar} functions of the CCB concentration $c_{\text{CCB}}$
\begin{equation}
    \label{eq:ccbs}
    \begin{aligned}
        \gamma_1(c_{\text{CCB}}) = \gamma_{1,\mathrm{max}} \left(1-p_1\dfrac{c_{\text{CCB}}^{2}}{c_{\text{CCB}}^{2}+c_{\mathrm{CCB},50}^{2}}\right), \\
        \gamma_3(c_{\text{CCB}}) = \gamma_{3,\mathrm{max}} \left(1-p_3\dfrac{c_{\text{CCB}}^{2}}{c_{\text{CCB}}^{2}+c_{\mathrm{CCB},50}^{2}}\right),
    \end{aligned}
\end{equation}
wherein $\gamma_{1, \max}$ and $\gamma_{3, \max}$ are the maximum attainable values, whereas $p_1$ and $p_3$ are parameters and $c_{\text{CCB},50}$ is the half activation CCB concentration.

\subsection{Smooth muscle cell  activation model}
From the work of Uhlmann and Balzani \cite{uhlmann2022simulation2} we have the active arterial response
\begin{equation}
    \Psi_{\mathrm{a}}^{(f)} =
    \frac{\mu_{\mathrm{a}}}{2}
    \left( n_{\mathrm{C}}^{(f)} + n_{\mathrm{D}}^{(f)} \right)
    (\lambda_{\mathrm{e}}^{(f)} - 1)^2,
\end{equation}
wherein, $\mu_{\mathrm{a}}$ is a stiffness constant, $n_{\mathrm{C}}^{(f)}$ and $n_{\mathrm{D}}^{(f)}$ are the fraction of myosin heads in states C and D, respectively. The stretch $\lambda^{(f)}$ is multiplicatively split into an elastic $\lambda_{e}^{(f)}$ and an active part $\lambda_{\mathrm{a}}^{(f)}$
\begin{equation}
    \lambda^{(f)} = \lambda_{\mathrm{e}}^{(f)} \lambda_{\mathrm{a}}^{(f)}.
\end{equation} 
The rate of change of active strain is given by
\begin{equation}
    \dot{\lambda}_{\mathrm{a}}^{(f)} =
    \beta_1 \dfrac{P_{\mathrm{a}}^{(f)} - P_{\mathrm{c}}^{(f)}}
    {P_{\mathrm{a}}^{(f)} - \beta_2},
\end{equation}
where $\beta_1$ and $\beta_2$ are material parameters, $P_{\mathrm{a}}^{(f)}$ and $P_{\mathrm{c}}^{(f)}$ are the active and driving stresses, respectively and are given by
\begin{equation}
    \begin{aligned}
        P_{\mathrm{a}}^{(f)} =
        \mu_{a} \left( n_{\mathrm{C}}^{(f)} + n_{\mathrm{D}}^{(f)} \right)
        (\lambda_{\mathrm{e}}^{(f)} - 1), \\
        P_{\mathrm{c}}^{(f)} =
        \kappa n_{\mathrm{C}}^{(f)},
    \end{aligned}
\end{equation}
where $\kappa$ is the maximum driving stress.

\section{Coupled diffusion-deformation problem} \label{sec:coupling}
The problem is described by the coupling of balance of linear momentum with a diffusion-reaction equation describing the local distribution of $n$ different drug concentrations $c_1 ,\dots , c_n$. 
The set of partial differential equations is given as: 
\begin{align}
&\text{Balance of linear momentum: } & \rho_s \frac{\partial ^{2} \mathbf u}{\partial t^{2}} -\nabla \cdot \mathbf P \, ( \mathbf u ,c_1,\dots,c_n) 
&= \rho_s \mathbf b  \label{eq:momentumBalanceStrong}
\\
&\text{Diffusion-reaction of drug i: } & \frac{\partial c_i}{\partial t} -\nabla \cdot \mathbf D \nabla c_i - r(c_i) 
&= 0 
\label{eq:reactionDiffusionStrong}
\end{align}
The function $r(c_i)$ describes, for example, the metabolism rate resulting in the reduction of drug~$i$. 
$\mathbf P $ is the first Piola-Kirchhoff stress tensor, $\rho_s$ the density, and $\boldsymbol{D}$ the diffusion tensor. 
Here, both processes occur simultaneously, namely the deformation due to blood pressure variations and the diffusion of drugs through the arterial wall. 

The set of equations is solved using the finite element method. We apply the Galerkin method to obtain the weak form of \cref{eq:momentumBalanceStrong,eq:reactionDiffusionStrong} with test functions $\mathbf v$ and $v$ for $\Omega \in \mathbb{R}^{3}$. For a more detailed derivation, we refer to \cite{wriggers2008nonlinear}. 
With $\mathbf  C_v = \nabla \mathbf  v^{T} \mathbf  F + \mathbf  F^{T} \nabla \mathbf  v $ and second Piola-Kirchhoff stress tensor $\mathbf  S$ we obtain the weak form of the balance of linear momentum
\begin{align}
\begin{split}
G^{\text{u}}(\mathbf  u,c,\mathbf  v) :=& \int_{\hat{\Omega}} \rho_s \frac{\partial^{2} \mathbf  u }{\partial t^{2}} \cdot \mathbf  v ~ \text{d} \hat{x} + \int_{\hat{\Omega}} \frac{1}{2} \mathbf  C_{\mathbf  v} :\mathbf  S  ~\text{d} \hat{x} -  \int_{\hat{\Omega}} \rho_s \mathbf  b \cdot \mathbf  v ~\text{d} \hat{x} - \int_{\partial \hat \Omega} \mathbf  v \cdot \mathbf  T_0 ~\text{d} \sigma =0.
\end{split}
\end{align}
Here, on the Neumann boundary $\mathbf T_0 = \mathbf{P} \cdot \mathbf n_0$ describes the stresses in the outward normal direction. We apply deformation-dependent loads as described in chapter 4.2.5 of \cite{wriggers2008nonlinear}, where the corresponding pressure values are inserted. Note, that this generates another nonlinear component to the tangent matrix, which is included in the structure block.

Transforming the volume elements to the reference configuration $\text{d} x = J \text{d} \hat{x}$ and mapping $  \mathbf F_R^{-1}: \mathbb{R}^{3 \times 3} \rightarrow \hat \Omega$, with $\nabla (\phi \circ \mathbf F_R^{-1})(x) = \nabla \phi(\mathbf F_R^{-1}(x)) \mathbf F^{-1}$ and deformation gradient $\mathbf F$, similarly, we obtain the weak form of the diffusion-reaction equation:
\begin{align}
\begin{split}
G^{\text{c}}(\mathbf u,c,v) :=& \int_{\hat{\Omega}} \frac{\partial c}{\partial t} \cdot  v~ \text{d} \hat{x} + \int_{\hat{\Omega}} \nabla c \cdot  \mathbf D_F \cdot \nabla v  - r(c)~ v ~J~\text{d} \hat{x} - \int_{\partial \hat{\Omega}}  \mathbf D_F  \cdot \nabla c \cdot \mathbf n_0~v ~J~\text{d} \hat{x}=0,
\end{split}
\end{align}
 \text{ with } $\mathbf D_F := \mathbf F^{-1} \mathbf D \mathbf F^{-T} $.
 
The coupling conditions reflect the influence of the drug on the deformation process and the extent of the deformation on the diffusion process. 
Generally, the diffusion process is influenced by the deformation via the change in element volume. More precisely, the determinant of the deformation gradient changes when the geometry deforms. For almost incompressible material the influence diminishes.
The influence of drug concentration is reflected depending on the underlying modeling of the arterial wall. We have nonlinear dependencies between the governing equations.
Thus, the equations are linearized using the Newton-Raphson scheme.

Then, we define the residual $\mathcal{R}$ with $G^{\text{u}}(\mathbf u, c, \mathbf v)$ and $G^{\text{c}}(\mathbf u, c ,v)$ for Newton iterations $k=0,1,\dots$ and the linearization $\Delta$ of the weak forms $G^{\text{c}}$ and $G^{\text{u}}$ in the Jacobian $\mathcal{J}$ respectively.
\begin{equation*}
\mathcal{R}(\mathbf u^{k},c^{k})
:= 
	\left.
	\begin{pmatrix} 
		G^{\text{u}}(\mathbf u, c, \mathbf v) \\
		G^{\text{c}}(\mathbf u, c ,v)
	\end{pmatrix}
	\right\rvert_{(\mathbf u^{k},c^{k})},
	\quad
	\mathcal{J}(\mathbf u^{k},c^{k})
:=
	\begin{pmatrix}
		\Delta_\text{u} [\mathbf G^{\text{u}}] & ~~ \Delta_\text{c} [\mathbf G^{\text{u}}] \\
		\Delta_\text{u} [\mathbf G^{\text{c}}] & ~~ \Delta_\text{c } [\mathbf G^{\text{c}}] 
	\end{pmatrix}
=
	\left.
	\begin{pmatrix}
		\mathbf K^{\text{u} \text{u}} & \mathbf K^{\text{u} \text{c}} \\
		\mathbf K^{\text{c} \text{u}} & \mathbf  K^{\text{c} \text{c}} 
	\end{pmatrix}
	\right\vert_{(\mathbf u^{k},c^{k})}.
\end{equation*}
Ultimately, the Newton solve for the increments $\delta \mathbf u^{k+1}$ and $\delta c^{k+1}$ is defined as:
\begin{equation}
\mathcal{J}(\mathbf u^{k},c^{k})
\begin{pmatrix}
\delta \mathbf u^{k+1} \\
\delta c^{k+1}
\end{pmatrix}
=
-\mathcal{R}(\mathbf u^{k},c^{k})
\label{eq:NewtonSolve}
\end{equation}

We can extend~\cref{eq:NewtonSolve} with respect to the time discretization. 
Time steps are labeled with $n$.
The residual components are included in $\boldsymbol{r}^{\text{u}}$ and $\boldsymbol{r}^{\text{c}}$. 
The following system of equations is implicit with respect to the coupling of diffusion and deformation:
\begin{align}
\begin{split}
\begin{bmatrix}
\boldsymbol{K}^{\text{uu}}_{n+1,k} &\boldsymbol{K}^{\text{uc}}_{n+1,k}\\
\boldsymbol{K}^{\text{cu}}_{n+1,k} & \boldsymbol{K}^{\text{cc}}_{n+1,k}
\end{bmatrix}
\begin{bmatrix}
\delta \mathbf u^{n+1,k+1} \\
\delta c^{n+1,k+1}
\end{bmatrix}
+
\begin{bmatrix}
\boldsymbol{M}_{n+1,k}^{\text{u}} & \bzero\\
\bzero           & \bzero
\end{bmatrix}
\begin{bmatrix}
{\delta\Ddot{\mathbf u}}^{n+1,k+1}\\
\bzero
\end{bmatrix} 
+
\begin{bmatrix}
\bzero & \bzero\\
\bzero & \boldsymbol{M}_{n+1,k}^{\text{c}}
\end{bmatrix}
\begin{bmatrix}
\bzero \\
\delta\dot c_{}^{n+1,k+1}
\end{bmatrix}
+ 
\begin{bmatrix}
\bbr_{n+1,k}^{\text{u}} \\
\bbr_{n+1,k}^{\text{c}}
\end{bmatrix}
   = \boldsymbol{0}.
\end{split}
\label{eq:linearizedSystemMatrix}
\end{align}

Finally, if we neglect the dynamic component and use a backward Euler scheme to discretize $c$ in time, we obtain the following fully implicit scheme from \eqref{eq:linearizedSystemMatrix}.
\begin{align}
\begin{array}{c}
		\left[
		\begin{array}{c l}
		\boldsymbol{K}^{\text{uu}}_{n+1,k} &\boldsymbol{K}^{\text{uc}}_{n+1,k}\\
		\boldsymbol{K}^{\text{cu}}_{n+1,k} & \boldsymbol{K}^{\text{cc}}_{n+1,k} - \frac{1}{\Delta t}  \boldsymbol{M}_{n+1,k}^{\text{c}}
		\end{array}
		\right]
		\begin{bmatrix}
\delta \mathbf u^{n+1,k+1} \\
\delta c^{n+1,k+1}
\end{bmatrix}
		+ 
		\left[
		\begin{array}{c}
		\bbr_{n+1,k}^{\text{u}} \\
		\bbr_{n+1,k}^{\text{c}} + \frac{1}{\Delta t}  \boldsymbol{M}_{n+1,k}^{\text{c}} \bd_{n}^{\text{c}} 
		\end{array}
		\right]
 		  = \boldsymbol{0}
  		\end{array}
    \label{eq:implicit discretized system}
    \end{align}
    
We use a $\mathcal{P}_2 - \mathcal{P}_2$ finite element discretization.

\section{Two-level overlapping Schwarz Preconditioners} \label{sec:two-level schwarz}

To efficiently solve the nonsymmetric system~\cref{eq:implicit discretized system}, we employ a preconditioned generalized minimal residual (GMRES) method~\cite{saad_gmres_1986}. 
The convergence of the (unpreconditioned) GMRES method will deteriorate when the mesh is refined. 
Hence, we will use preconditioning techniques to improve the convergence and to obtain parallel scalability for large problems. 
In particular, we use two-level overlapping Schwarz domain decomposition preconditioners with generalized Dryja--Smith--Widlund (GDSW) coarse spaces \cite{dohrmann_domain_2008,dohrmann_design_2017}. 
We use the FROSch (Fast and Robust Overlapping Schwarz)~\cite{Heinlein:2018:FRO} domain decomposition solver package, which is part of the Trilinos~\cite{trilinos,trilinos-website} library, for the implementation of these preconditioners; cf.~\Cref{sec:frosch}.

For the sake of simplicity, let us first describe the two-level Schwarz preconditioners just for a single-physics problem, for instance, for a diffusion problem 
$$
	\boldsymbol{K} \boldsymbol{c} = \boldsymbol{f}
$$
resulting from a finite element discretization with a finite element space $V$. 
In~\Cref{sec:monolithic Schwarz preconditioners}, we discuss how to construct the preconditioners to block systems of the form~\cref{eq:implicit discretized system}. 
The preconditioners are then based on an overlapping domain decomposition of the computational domain $\Omega$ into $N$~overlapping subdomains $\left\lbrace \Omega_i' \right\rbrace_{i=1,\ldots,N}$. 
In practice, we construct the overlapping domain decomposition from a nonoverlapping domain decomposition $\left\lbrace \Omega_i \right\rbrace_{i=1,\ldots,N}$ by recursively extending the subdomains by layers of finite elements; 
when extending the subdomains by $k$~layers of elements with a width of $h$, the overlap has a width of $\delta = kh$. 
The overlapping subdomains then conform with the computational mesh on $\Omega$. 
Note that in our implementation (cf.~\Cref{sec:feddlib trilinos}), we construct the overlap algebraically, via the graph induced by the nonzero pattern of~$K$. 
We will denote the size of the algebraically determined overlap with $\widehat{\delta}=k$, for $k$~layers of degrees of freedom. 
The initial nonoverlapping domain decomposition can be obtained using mesh partitioning algorithms, for instance, using the parallel ParMETIS library~\cite{karypis1997parmetis}. 

Let $V_i$ be the finite element space on the local overlapping subdomain $\Omega_i'$ with homogeneous zero boundary conditions on $\partial\Omega_i'$, and let $R_i \colon V \to V_i$ be the restriction from the global to the local finite element space on $\Omega_i'$. 
The corresponding prolongation operators are then given by $R_i^\top$. 
Furthermore, we formally introduce a (global) coarse space $V_0$ spanned by coarse basis functions corresponding to the columns of the matrix $\Phi \colon V_0 \to V$; 
we will introduce specific coarse spaces in detail in~\Cref{sec:coarse spaces}. 

The two-level additive overlapping Schwarz preconditioner reads
\begin{equation} \label{eq:two-level schwarz}
	M^{-1} = \underbrace{\Phi K_0^{-1} \Phi^\top}_{\text{second level}} + \underbrace{\sum\nolimits_{i=1}^{N} R_i^\top K_i^{-1} R_i}_{\text{first level}},
\end{equation}
where $K_0$ is the coarse matrix, 
\begin{align}\label{eq:coarseMatrix}
K_0 = \Phi^\top K \Phi,
\end{align}
and  $K_i = R_i K R_i^\top$ are the local matrices. 
Note that the coarse level corresponds to solving a globally coupled coarse problem, whereas the first level corresponds to solving $N$~decoupled local problems that can be solved concurrently in parallel computations. 
For more details on standard Schwarz preconditioners, we refer to~\cite{smith_domain_1996,toselli_domain_2005}.

\subsection{GDSW and RGDSW coarse spaces for overlapping Schwarz}\label{sec:coarse spaces}
\label{subsection:GDSWandRGDSW}

The coarse level, which is determined by the choice of the coarse space, is essential for the numerical scalability of the solver. In particular, without a suitable coarse level, the convergence of the Krylov method with a one-level Schwarz preconditioner will increase with an increasing number of subdomains. Here, we will focus on coarse spaces that can be constructed from the parallel distributed system matrix, requiring little to no additional information.

In particular, we consider GDSW~\cite{dohrmann_family_2008,dohrmann_domain_2008} and reduced GDSW (RGDSW)~\cite{dohrmann_design_2017} coarse spaces. 
Let us discuss both approaches in a common algorithmic framework. 
Both approaches start with a decomposition into nonoverlapping subdomains $\left\lbrace \Omega_i \right\rbrace_{i=1,\ldots,N}$, as already mentioned in~\Cref{sec:two-level schwarz}. 
First, we define the interface of the nonoverlapping domain decomposition as
$$
	\Gamma = \bigcup_{i \neq j} \left( \partial\Omega_i \cap \partial\Omega_j \right) \setminus \partial\Omega_D,
$$
where $\partial\Omega_D$ is the global Dirichlet boundary. 
Then, we decompose the interface $\Gamma$ into components $\gamma_k \subset \Gamma$ such that
\begin{equation} \label{eq:conver_interface}
	\Gamma = \bigcup_k \gamma_k.
\end{equation}
This decomposition may be disjoint ($\gamma_i\cap\gamma_k=\emptyset$ for $i\neq k$) or overlapping. 
A typical choice is to decompose the interface into faces, edges, and vertices, where (in three dimensions)
\begin{itemize}
	\item a face is a set of nodes that belong to the same two subdomains,
	\item an edge is a set of at least two nodes that belong to the same (more than two) subdomains,
	\item a vertex is a single node that belongs to the same (more than two) subdomains.
\end{itemize}
Then, we define functions $\phi_l^{\gamma_k}$ on $\Gamma$ with $\text{supp} (\phi_l^{\gamma_k}) \subset \gamma_k$ such that
$$
	N_\Gamma \subset \text{span} \left\lbrace \phi_l^{\gamma_k} \right\rbrace_{l,k},
$$
where $N_\Gamma$ is the restriction of the null space~$N$ of the global Neumann matrix to the interface $\Gamma$. 
The global Neumann matrix is obtained by assembling the original problem on~$\Omega$ but with $\partial\Omega_D=\emptyset$, i.e., a Neumann condition on~$\partial\Omega$. 
Then, we gather the functions as columns of a matrix
$$
	\Phi_\Gamma = \begin{pmatrix} \phi_1^{\gamma_1} & \phi_2^{\gamma_1} & \cdots & \phi_1^{\gamma_2} & \cdots \end{pmatrix}.
$$ 
The matrix $\Phi$ is obtained from computing the discrete energy-minimizing extension of the interface values $\Phi_\Gamma$ to the interior of the subdomains
\begin{equation} \label{eq:Phi}
	\Phi 
	=
	\begin{bmatrix}
		\Phi_I \\
		\Phi_\Gamma
	\end{bmatrix}
	=
	\begin{bmatrix}
		- \boldsymbol{K}_{II}^{-1} \boldsymbol{K}_{I \Gamma} \Phi_\Gamma \\
		\Phi_\Gamma
	\end{bmatrix},
\end{equation}
where the index~$I$ indicates all degrees of freedom that do not correspond to the interface $\Gamma$, and $\boldsymbol{K}_{II}$ and $\boldsymbol{K}_{I \Gamma}$ are obtained by reordering and partitioning $\boldsymbol{K}$ as follows:
$$
	\boldsymbol{K}
	=
	\begin{bmatrix}
		\boldsymbol{K}_{II} & \boldsymbol{K}_{I\Gamma} \\
		\boldsymbol{K}_{\Gamma I} & \boldsymbol{K}_{\Gamma\Gamma}
	\end{bmatrix}.
$$

As mentioned before, we will employ GDSW and RGDSW coarse spaces in this paper. 
Based on the algorithmic framework described above, they can be easily introduced. 

In the GDSW coarse space, we choose the interface components $\gamma_k^{\text{GDSW}}$ as faces, edges, and vertices. 
Then, the functions $\phi_l^{\gamma_k^{\text{GDSW}}}$ are defined as the restrictions of a basis of the null space~$N$ to the interface components $\gamma_k^{\text{GDSW}}$. 
For example, in the case of a scalar diffusion problem, the null space consists only of constant functions, and we can set $\phi_l^{\gamma_k^{\text{GDSW}}}$, $l=1$, to one on $\gamma_k^{\text{GDSW}}$ and zero elsewhere. 

For a homogeneous two-dimensional diffusion problem with irregular subdomains as, for instance, obtained by METIS, the two-level Schwarz preconditioner~\cref{eq:two-level schwarz} with the GDSW coarse space yields the condition number estimate
$$
	\kappa \left( \boldsymbol{M}^{-1} \boldsymbol{K} \right)
	\leq
	C \left(1 + \frac{H}{\delta}\right) \left( 1 + \log \left( \frac{H}{h} \right) \right)^2;
$$
cf.~\cite{dohrmann_domain_2008}. Here, $h$ is the typical diameter of the finite elements, $H$ the diameter of the subdomains, and $\delta$ is the overlap.

The construction of the RGDSW coarse space employed in this paper (``Option 1'' in~\cite{dohrmann_design_2017}) is based on a different, overlapping decomposition of the interface, generally resulting in a coarse space of lower dimension. 
In particular, let us consider the interface components $\gamma_k^{\text{GDSW}}$ from the GDSW coarse space, which correspond to faces, edges, and vertices, and let $\mathcal{S}^{\gamma_k^{\text{GDSW}}}$ be the set of subdomains that $\gamma_k^{\text{GDSW}}$ belongs to. 
Then, we select all interface components $\gamma_k^{\text{GDSW}}$ with
$$
	\nexists \gamma_l^{\text{GDSW}}, l\neq k: \mathcal{S}^{\gamma_k^{\text{GDSW}}} \subset \mathcal{S}^{\gamma_l^{\text{GDSW}}},
$$
and we denote them by $\widehat{\gamma}_k^{\text{RGDSW}}$. 
To satisfy~\cref{eq:conver_interface}, we choose
$$
	\gamma_k^{\text{RGDSW}} = \widehat{\gamma}_k^{\text{RGDSW}} \cup \bigcup_{\mathcal{S}^{\gamma_l^{\text{GDSW}}} \subset \mathcal{S}^{\widehat{\gamma}_k^{\text{RGDSW}}}} \gamma_l^{\text{GDSW}},
$$
that is, we add all edges, faces, and vertices that belong to a subset of the subdomains that $\widehat{\gamma}_k^{\text{RGDSW}}$ belongs to. 
To define the corresponding interface functions $\phi_l^{\gamma_k^{\text{GDSW}}}$, we define a scaling function $s\colon \Gamma \to \mathbb{R}$ with 
$$
	s(n) = \left| \left\lbrace \gamma_k^{\text{RGDSW}} \vert n \in \gamma_k^{\text{RGDSW}}  \right\rbrace \right|
$$
for $n \in \Gamma$. 
Denoting by $\widehat{\phi}_l^{\gamma_k^{\text{RGDSW}}}$ restrictions of a basis of the null space $N$ to the interface components $\gamma_k^{\text{RGDSW}}$, we define
$$
	\phi_l^{\gamma_k^{\text{RGDSW}}} (n) = \frac{\widehat{\phi}_l^{\gamma_k^{\text{RGDSW}}} (n)}{s(n)}
$$
for $n \in \Gamma$; 
this means that in each node~$n$ we scale the restrictions of the basis functions of the null space by the inverse of the number of RGDSW interface components the node belongs to. For a simple homogeneous diffusion problem, the two-level Schwarz preconditioner in~\cref{eq:two-level schwarz} with the RGDSW coarse space yields the condition number estimate
$$
\kappa \left( \boldsymbol{M}^{-1} \boldsymbol{K} \right)
\leq
C \left(1 + \frac{H}{\delta}\right) \left( 1 + \log \left( \frac{H}{h} \right) \right);
$$
cf.~\cite{dohrmann_design_2017}.

\subsection{Monolithic overlapping Schwarz preconditioners}\label{sec:monolithic Schwarz preconditioners}

For preconditioning block systems of the form~\cref{eq:implicit discretized system}, we employ monolithic Schwarz preconditioning techniques; cf.~\cite{klawonn_overlapping_1998,klawonn_comparison_2000}. In particular, we use monolithic GDSW preconditioners
\begin{equation} \label{eq:monolithic:two-level schwarz}
	\mathcal{M}^{-1} 
	= 
	\phi \mathcal{K}_0^{-1} \phi^\top + \sum\nolimits_{i=1}^{N} \mathcal{R}_i^\top \mathcal{K}_i^{-1} \mathcal{R}_i,
\end{equation}
as introduced in~\cite{heinlein_monolithic_2019,heinlein_reduced_2020} for a block matrix
\begin{equation} \label{eq:mcK}
	\mathcal{K}
	=
	\begin{bmatrix}
		\mathcal{K}_{u,u} & \mathcal{K}_{u,c} \\
		\mathcal{K}_{c,u} & \mathcal{K}_{c,c}
	\end{bmatrix}.
\end{equation}
In~\cref{eq:monolithic:two-level schwarz}, the local subdomain matrices are given by $\mathcal{K}_i = \mathcal{R}_i \mathcal{K} \mathcal{R}_i^\top$ with restriction operators
$$
	\mathcal{R}_i
	=
	\begin{bmatrix}
		\mathcal{R}_{u,i} & 0 \\
		0 & \mathcal{R}_{c,i}
	\end{bmatrix},
$$
where $\mathcal{R}_{u,i}$ and $\mathcal{R}_{c,i}$ correspond to the restriction operators of the $u$ and $c$ degrees of freedom, respectively, for $\Omega_i'$. 
Also the coarse matrix $\mathcal{K}_0 = \phi^\top \mathcal{K} \phi,$ has a block structure given by the matrix 
$$
	\phi
	=
	\begin{bmatrix}
		\phi_{u,u_0} & \phi_{u,c_0} \\
		\phi_{c,u_0} & \phi_{c,c_0}
	\end{bmatrix},
$$
which has the coarse basis functions as its columns.

It is constructed analogously to~\cref{eq:Phi} based on the null spaces of the Neumann matrices corresponding to the diagonal blocks of $\mathcal{K}$; cf.~\cref{eq:mcK}. See also~\cite{heinlein_frosch_2022} for monolithic Schwarz preconditioners for land ice problems, which have a similar matrix structure as~\cref{eq:implicit discretized system}.

\subsection{Recycling strategies}\label{sec:Recycling Strategies}

To save computing , information from previous time steps or Newton iterations can be reused. 
Specifically, three recycling strategies will be analyzed further: 
the reuse of the symbolic factorizations of local subdomain matrices~$\mathcal{K}_i$ and of the coarse matrix~$\mathcal{K}_0$---we denote this strategy by (SF)---the reuse of the coarse basis~$\Phi$, which is used to assemble~$\mathcal{K}_0$ (cf.~\cref{eq:coarseMatrix})---we denote this strategy by (CB)---and the reuse of the entire coarse matrix~$K_0$, denoted by (CM).

Generally, the reuse of a symbolic factorization is advantageous as the nonzero pattern of the system typically stays the same. 
Reusing the entire coarse matrix~$\mathcal{K}_0$ not only eliminates the time for setting it up but also the time required for its factorization. 
The reuse of the coarse basis and the coarse matrix, respectively, are more invasive strategies. 
They can further decrease the setup time but can also increase the iteration count, as the preconditioner no longer adapts exactly to the new Newton iteration or time step. Note, however, that we are still solving the correct tangent matrix system in each Newton step since only the preconditioner is affected.

\section{Software Ecosystem}

For the large-scale simulations of arterial walls that are characterized by the models in \Cref{sec:modeling}, we use the different software libraries AceGen, Trilinos, FEDDLib, and FROSch together with a newly developed AceGen--FEDDLib interface. 
As the main software package, we use the library FEDDLib~\cite{feddlib} (\textbf{F}inite \textbf{E}lement and \textbf{D}omain \textbf{D}ecomposition \textbf{Lib}rary), and for some functionalities, we make use of other libraries. In particular, we employ data services, parallel linear algebra, solvers, and preconditioners to solve our systems efficiently from the open-source software library Trilinos~\cite{trilinos}. For the implementation of the solid material models, we use the commercial symbolic software AceGen, which is based on Mathematica. More specifically, we generate the code that implements the material models using AceGen and call it through an interface between AceGen and FEDDLib.

Through an interface between AceGen and FEDDLib, which has been developed within the present project, the assembled structure-chemical interaction system can be passed along. 
FEDDLib then uses the solvers and preconditioners provided by Trilinos to solve the system. 
We especially focus on the Trilinos package FROSch (\textbf{F}ast and \textbf{R}obust \textbf{O}verlapping \textbf{Sch}warz), which includes the parallel implementation of the two-level overlapping Schwarz preconditioner as described in \Cref{sec:two-level schwarz}.

\subsection{FEDDLib \& Trilinos} \label{sec:feddlib trilinos}

FEDDLib~\cite{feddlib} is a \Cpp-based, object-oriented finite element library built on top of the Trilinos software infrastructure. 
Trilinos~\cite{trilinos,trilinos-website} is a software framework containing mathematical software libraries for the solution of large-scale, complex multiphysics engineering and scientific problems; cf.~\cite{heroux_overview_2005}. It is organized as a collection of (mostly) interoperable packages. The packages are categorized in the six product areas \emph{Data Services}, \emph{Linear Solvers}, \emph{Nonlinear Solvers}, \emph{Discretizations}, \emph{Framework}, and \emph{Product Manager}. 

The main design principle of FEDDLib is to reach outstanding parallel efficiency and robustness for complex single and multiphysics applications by closely integrating finite element discretization and domain decomposition-based solvers. In particular, the finite element assembly in FEDDLib provides the data structures required to make use of the full potential of state-of-the-art domain decomposition algorithms; in this context, the main focus is on overlapping Schwarz solvers in Trilinos' domain decomposition solver package FROSch~\cite{Heinlein:2018:FRO}.

The main components of FEDDLib can be categorized into \emph{linear algebra}, \emph{mesh handling}, \emph{finite element assembly}, \emph{boundary condition handling}, \emph{linear} and \emph{nonlinear solvers}, \emph{time-stepping},  \emph{parallel I/O}, and \emph{application-specific classes}; an overview over the functionality in the year 2020 is given in the PhD thesis of C.~Hochmuth~\cite{kups11345}. Afterwards, further significant developments have been carried out, for instance, on adaptive mesh refinement by one of the authors~\cite{Sassmannshausen:2021:AFE} and an interface for AceGen generated code, which is briefly discussed in~\Cref{sec:acegen,sec:interface}. Application problems currently implemented in FEDDLib range from simple scalar elliptic problems, such as diffusion, to complex multiphysics problems, such as fluid-structure interaction with nonlinear models for the fluid and solid. 

FEDDLib provides wrappers and interfaces to many packages of Trilinos. First, the parallel linear algebra in FEDDLib is done via wrappers for the parallel linear algebra classes in Trilinos (product area \emph{Data Services}); in particular, FEDDLib uses Xpetra, which in turn is a lightweight interface for both the Epetra and Tpetra linear algebra frameworks in Trilinos. Whereas Epetra is the older linear algebra framework from the origins of Trilinos, Tpetra is the current linear algebra framework. One of the main advancements of Tpetra compared with Epetra is the integration of the performance portability libraries Kokkos~\cite{9485033,9502936,CarterEdwards20143202} and Kokkos Kernels~\cite{rajamanickam_kokkos_2021}. This also enables the efficient use of GPUs; see, for instance,~\cite{yamazaki_experimental_2023} for the application of FROSch on GPUs via Kokkos and Kokkos Kernels. 

Therefore, we mostly use the newer Tpetra package, and we will focus on the software stack based on Tpetra in our discussion here. In particular, FEDDLib provides interfaces to Trilinos' nonlinear solver package NOX (product area \emph{Nonlinear Solvers}) as well as several packages for linear solvers and preconditioners (product area \emph{Linear Solvers}), which we access via the unified solver interface Stratimikos~\cite{bartlett_stratimikos_2006}, which in turn makes use of the interoperability framework Thyra~\cite{bartlett_thyra_2007}. In particular, we employ the Belos and Amesos2~\cite{bavier_amesos2_2012} packages for iterative and direct solvers, respectively; note that the direct solvers---except for the serial direct solver KLU---are not part of Trilinos and can therefore only be called via the interface Amesos2. For preconditioning, we mostly use the domain decomposition package FROSch~\cite{Heinlein:2018:FRO} and the block-preconditioning package Teko~\cite{cyr_teko_2016}; other available options are Ifpack2, which combines one-level Schwarz methods with inexact factorization preconditioning techniques, and the algebraic multigrid package MueLu~\cite{berger-vergiat_muelu_2019}. In~\Cref{sec:frosch}, we will discuss the FROSch package in more detail since it implements the Schwarz domain decomposition preconditioners discussed in~\Cref{sec:two-level schwarz} and investigated in our~\Cref{sec:numerical results} on numerical results. Finally, the toolbox package Teuchos (product area \emph{Framework}), which provides smart pointers, parameter lists, and XML parsers, as well as the graph partitioning and load-balancing package Zoltan2 (product area \emph{Data Services}) for mesh partitioning are employed.

\subsection{FROSch} \label{sec:frosch}

The Trilinos package FROSch (Fast and Robust Overlapping Schwarz)~\cite{Heinlein:2018:FRO} provides a highly-scalable, parallel framework for overlapping Schwarz domain decompositions solvers; the first version of the implementation was described in~\cite{heinlein_parallel_2016a,heinlein_parallel_2016b} and the extension to monolithic GDSW and RGDSW preconditioners in~\cite{heinlein_monolithic_2019}.

The design of the FROSch package is based on the concept of combining Schwarz operators in additive or multiplicative ways, resulting in different types of Schwarz preconditioners; cf.~\cite{toselli_domain_2005,heinlein_parallel_2016a,Heinlein:2018:FRO}. The Schwarz preconditioners are constructed algebraically, that is, based on the sparsity pattern of the parallel-distributed system matrix. For the first level, the overlapping subdomains are constructed from nonoverlapping subdomains by extending the subdomains by layers of adjacent degrees of freedom. Adjacency is defined based on nonzero off-diagonal entries in the system matrix. 

For the coarse level, we employ extension-based coarse spaces, such as GDSW and RGDSW; cf.~\Cref{sec:coarse spaces}. 
The construction described in~\Cref{sec:coarse spaces} is algebraic except for the fact that the null space of the Neumann matrix is required. In particular, for elasticity problems, the null space is spanned by the rigid body modes, and the (linearized) rotations cannot be constructed algebraically; they have to be either provided by the user or constructed using coordinates of the finite element nodes. For more details on the algebraic construction of GDSW-type coarse spaces for elasticity problems, see~\cite{heinlein_fully_2021}.

FROSch preconditioners can be easily constructed via Trilinos' unified solver interface Stratimikos using only the parallel-distributed system matrix and a list of parameters defining the specific setup of the preconditioners. FROSch is based on the Xpetra wrappers, which wrap the Epetra and Tpetra linear algebra frameworks of Trilinos; cf.\ the discussion in~\Cref{sec:feddlib trilinos}. If Tpetra is used, the performance portability libraries Kokkos and Kokkos Kernels can be employed, enabling the use on, for instance, GPU architecture~\cite{yamazaki_experimental_2023}. Due to its algebraic implementation, FROSch preconditioners can be constructed recursively, resulting in multi-level Schwarz preconditioners; see~\cite{heinlein:2023:pst}.

\subsection{AceGen} \label{sec:acegen} 
AceGen \cite{korelc2016automation} is a Mathematica package that is used to automatically derive mathematical expressions required in numerical procedures.
AceGen exploits the symbolic and algebraic capabilities of Mathematica by combining them with a hybrid symbolic-automatic differentiation technique. 
AceGen uses automatic code generation and simultaneous optimization of expressions. 
Furthermore, AceGen's ability to generate code for multiple languages can be leveraged to generate finite element code to run simulations with, e.g., FEAP (Finite Element Analysis Program), Elfen, Abaqus, and the Mathematica-based finite element environment AceFEM. 
Using the AceGen--AceFEM combination, new finite elements can be rapidly developed and tested. 
Here, we use AceGen to generate finite element code with consistent linearizations, which is important for the optimal performance of the Newton algorithm.

\subsection{Interface} \label{sec:interface}
The newly developed AceGen interface enables the use of AceGen-generated AceFEM (AceGen--AceFEM) finite elements in external software libraries, by providing \Cpp~finite element classes that can be used to compute element-level quantities such as residuals, stiffness matrices, history variables and post-processing quantities.
The interface contains a set of functions that handle the manipulation of data structures that are required by the AceGen--AceFEM finite elements. It greatly simplifies the use of the generated finite elements by providing a unified object-oriented interface and obscures low-level memory management, which tends to be error-prone.
The interface theoretically also enables the use of the various finite elements available in the AceShare finite element library (which contains numerous community-contributed finite elements), in \Cpp-based finite element codes.

Accordingly, the corresponding FEDDLib assembly routines have been modified to include the externally generated AceGen--AceFEM finite elements through the AceGen interface, enabling their use in high-performance computing environments. An object-oriented, factory-based concept enables specific assembly routines to be built on top of a generalized base class, which defines the interface. Depending on the specific problem at hand, the base class' pure virtual functions, most importantly \texttt{getRhs()} and \texttt{getJacobian()}, are overridden in the derived class, with the corresponding finite element assembly routines. Mainly, FEDDLib only interacts with the base class to retrieve the element's residual vector (\texttt{getRhs()}) and Jacobian matrix (\texttt{getJacobian()}), that are necessary to setup the linearized Newton system as depicted in~\cref{eq:NewtonSolve}. 
The elementwise system matrix and vector---for the structure-chemical interaction, the Jacobian would be the 2$\times$2 block system---is then assembled into the parallely distributed \texttt{Tpetra::CrsMatrix} and \texttt{Tpetra::MultiVector} data types.
These are in turn passed along to the linear solvers and preconditioners in Trilinos.

\section{Numerical Results} \label{sec:numerical results}
For the numerical analysis of the coupling of the equations of linear momentum and reaction-diffusion (\Cref{sec:coupling}), we primarily focus on two aspects: the pharmaco-mechanical effects induced by the presence of medication and the weak and strong parallel scalability of the solvers. 

As described in \Cref{sec:coupling} for time discretization we apply the backward Euler scheme. As we neglect the dynamic component in the structural part, we employ different loading protocols for our numerical tests. Here, the load steps coincide with the time steps. 
For the discretization of \cref{eq:implicit discretized system}, we use tetrahedra with piecewise polynomials of degree two for both the displacement and the concentration.

Except for the weak scaling results in \Cref{sec:weak scaling}, all meshes are unstructured. 
We use Gmsh~\cite{geuzaine:2009:gmsh} to generate the unstructured meshes, and we partition them with METIS. 
The parallel results were obtained on the Fritz supercomputer at Friedrich-Alexander-Universität Erlangen-Nürnberg. 
It has 992~compute nodes, each with two Intel Xeon Platinum 8360Y ``Ice Lake'' processors (36~cores per chip, 2.4~GHz base frequency) and 256~GB of DDR4 RAM per node.

We use the Trilinos-based Newton solver \texttt{NOX} as the nonlinear solver (cf. \Cref{sec:feddlib trilinos}), which offers, among other features, different globalization and forcing term strategies for inexact Newton. 
If $J_{F_k} d_k = -F_k$ denotes the linear (nonsymmetric) system to be solved (see \Cref{eq:NewtonSolve}) in the $k$th Newton step, we use GMRES~\cite{saad_gmres_1986} (without restart) to find an approximate solution $\tilde{d}_k$ such that 
\begin{align}
\Vert J_{F_k} \tilde{d}_k + F_k \|_{l^2} \leq \eta_k \|F_k\|_{l^2}, 
\label{eq:Stopping Criterion GMRES}
\end{align}
where the forcing term (or tolerance) $\eta_k=10^{-6}$ is not changed between Newton steps. 
We employ no globalization strategies. 
In our implementation, this corresponds to restricting the line search in \texttt{NOX} to a full Newton step. 
Newton's method is stopped if the $l_2$-norm of the update $\tilde{d}_k$ is smaller than $10^{-8}$. 
The initial vector for the Newton iteration is the solution from the previous time step; 
the initial vector for GMRES is zero. 
The GMRES method is preconditioned with the two-level additive overlapping Schwarz method (\Cref{sec:two-level schwarz}), where the coarse space is either GDSW or RGDSW. 
As direct solvers to construct the preconditioner, we use Intel MKL Pardiso~\cite{schenk:2004:pardiso} on the subdomains and SuperLU\_DIST~\cite{superludist:2003} on the coarse level. 
SuperLU\_DIST is well-suited for the coarse solve, as we can execute it in parallel. Here, we use 5 processor cores to solve the coarse problem. 
Intel MKL Pardiso is used without threading. 
For the numerical results of this paper, we have always chosen the algebraically determined overlap $\widehat{\delta}$ as one layer of degrees of freedom.

\subsection{Pharmaco-mechanical results} \label{sec:Pharmaco-mechanicalResults}

\begin{figure}[tb] 
 \begin{subfigure}[b]{0.45\textwidth}
 	 	\centering
		\includegraphics[scale=2.3]{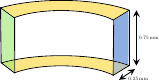}
	 	\centering	 	\includegraphics[scale=2.3]{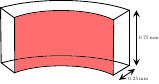}	
	 	\vspace{1cm}
	 \end{subfigure}
\begin{subfigure}[b]{0.49\textwidth}
	 	\centering
	 	\includegraphics[scale=0.95]{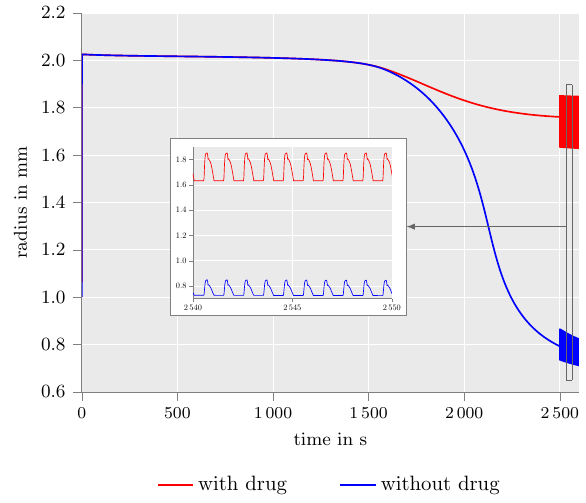}
	 \end{subfigure}
 	
	\caption{%
		Left: Arterial segment (quarter of tube) with inner radius 1\,mm, wall thickness 0.25\,mm and length 0.75\,mm. 
		Dirichlet boundary conditions: fixed in $x$ direction on $x=1$ plane (green), fixed in $y$ direction on $y=1$ plane (blue), fixed in $z$ direction on $z=0$ and $z=0.75$ planes (yellow). 
		On the interior wall (red), pressure is applied. 
		Right: Comparison of the response of the arterial radius with and without drug diffusion, with a	detailed extract for the time between 2\,540\,s and 2\,550\,s. 
		Our material model and the protocol in \Cref{fig:Loading Protocol long term} are used; 
		the drug diffusion is turned off for one of the simulations. 
		The oscillations in the vessel radius from 2\,500\,s onwards are caused by the onset of heartbeats.
	\label{fig:Arterial Wall and PharmacoMechanical Results}
	}
\end{figure}

\begin{figure}[tb]
	\centering
	\includegraphics[scale=0.9]{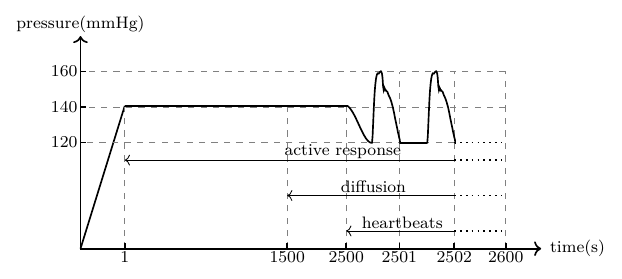}
	\caption{Loading protocol for long-term calculations with active response, diffusion process, and 100 heartbeats. Loading phase from $t=0$\,s to $t=1$\,s with 20 load steps. A constant pressure of 140\,mmHg is applied until $t=2\,500$\,s with $\Delta t=0.25$\,s. Active response starts at $t=1$\,s. Drug diffusion starts at $t=1\,500$\,s. Starting at $t=2\,500$\,s, a total of 100~heartbeats are simulated with $\Delta t=0.02$\,s. The simulation ends at $t=2\,600$\,s. The time-dependent change of the loading during the heartbeats was generated based on the flow rate profile from~\cite{balzani:2015:nmf}, which is based on data from~\cite{hemolab}.} 
	\label{fig:Loading Protocol long term}
\end{figure}

For the analysis of pharmaco-mechanical effects, we consider an idealized arterial segment with inner radius 1\,mm, outer radius 1.25\,mm and axial length of 0.75\,mm; see~\Cref{fig:Arterial Wall and PharmacoMechanical Results} (left). 
The two fiber directions are crosswise helically wound and lie at $30^{\circ}$ to the circumferential axis on the axial-circumferential plane, at every material point. 
The simulation protocol is given in \Cref{fig:Loading Protocol long term}. 

Between 0\,s and 1\,s, we have a linear increase in pressure to 140\,mmHg, which is then held constant until 2\,500\,s. 
Since simulating heartbeats that have systolic and diastolic pressures ranging from 120\,mmHg to 160\,mmHg from the beginning would be computationally expensive, we use a constant average pressure of 140\,mmHg instead. 
Starting from 1\,s, the active response of the material is turned on, and the evolution equations stated in \Cref{sec:modeling} are solved using the backward Euler scheme at every time step. 
It takes about 2\,000\,s to 2\,300\,s for the various internal variables to reach a steady state, which may be understood as corresponding to the physiological range. 
We observe an initial slow contraction, followed by a rapid contraction phase, as the fraction of myosin heads in state C increases. 
Drug diffusion begins at 1\,500\,s, signifying calcium channel blocker intake. 
We choose a diffusion coefficient of $D=6 \cdot 10^{-5}$. 
The Dirichlet boundary condition emulating drug inflow is implemented by switching the concentration from zero to one at $t=1\,500$\,s.
Starting from 2\,500\,s, at which point there is mature but incomplete diffusion of drugs, 100~heartbeats are simulated. 

The material parameters for the fully active material response are obtained from Uhlmann and Balzani~\cite{uhlmann2023chemo}. 
These material parameters were fitted by the authors to experimental results from Johnson \etal~\cite{johnson2009ca2+}, who investigated the contraction of a middle cerebral rat artery by applying a sequence of intravascular pressure with increasing pressure values. 
The material parameters of the passive response are listed in \Cref{tab:passive}, whereas the material parameters for the active part of the material model are given in \Cref{tab:chem_nOpt,tab:opt}. 

\begin{table}[tb]
	\caption{Passive material parameters \label{tab:passive}}
	\begin{center}
	\begin{tabular*}{\textwidth}{@{\extracolsep\fill}lccccc@{\extracolsep\fill}}
		\toprule 
	Parameter & $\alpha_1$ & $\alpha_2$ & $\alpha_3$ & $\alpha_4$ & $\alpha_5$ \\ \hline
	Value & $11.52507$\,kPa & $151.73775$\,kPa & $2.75662$ & $1.27631$\,kPa & $3.08798$ \\ \bottomrule
	\end{tabular*}
	\end{center}
\end{table}

\begin{table}[tb]
	\caption{Parameters for the chemical model \label{tab:chem_nOpt}}
		\begin{tabular*}{\textwidth}{@{\extracolsep\fill}lccccccccc@{\extracolsep\fill}}
		\toprule
	Parameter & $k_3$ & $k_4$ & $k_7$ & $Ca_{50}$ & $\gamma_2$ & $\gamma_3$ & $\lambda_{50, \, \text{c}}$ & $\bar{\lambda}_{\text{c}, \, \text{start}}$ & $\lambda_{\text{a}, \, \text{start}}$ \\ \hline
	Value & $0.134\,\text{s}^{-1}$ & $0.00166\,\text{s}^{-1}$ & $0.000066\,\text{s}^{-1}$ & $0.4\,\mu\text{M}$& $50\,\mu\text{M}^{-1}$ & $0.9\,\mu\text{M}$ & $1.2$ & $1.0$ & $1.0$\\ \toprule
	Parameter & $\gamma_4$ & $\zeta_1$ & $\gamma_5$ & $\zeta_2$ & $\Delta \bar{\lambda}_{\text{p}, \, \text{min}}$ & $\gamma_6$ & $\lambda_{50, \, \text{p}}$ & $\bar{\lambda}_{\text{p}, \, \text{start}}$ & $\lambda_{\text{e}, \, \text{start}}$\\ \hline
	Value & $200$ & $100\,\text{s}$ & $50\,\text{s}$ & $1000$ & $-0.00001$ & $1.5\,\text{s}^{-1}$ & $1.0$ & $1.0$ & $1.0$ \\ \bottomrule
	\end{tabular*}
\end{table}

\begin{table}[tb]
	\caption{Parameters for the active chemo-mechanical material model\label{tab:opt}}
	\begin{tabular*}{\textwidth}{@{\extracolsep\fill}lcccccc@{\extracolsep\fill}}
\toprule
	Parameter &$\eta$ & $\gamma_1$ & $\dot{\bar{\lambda}}_{\text{c}, \, \text{max}}$ & $\dot{\bar{\lambda}}_{\text{c}, \, \text{min}}$ & $\dot{k}_{2/5, \, \text{max}}$ & $\dot{k}_{2/5, \, \text{min}}$ \\ \hline
	Value & $0.1624\,\text{s}^{-1}$ & $0.5131\,\mu\text{M}$ & $0.0443\,\text{s}^{-1}$ & $-0.0443\,\text{s}^{-1}$ & $0.0009735\,\text{s}^{-2}$ & $-0.0010694\,\text{s}^{-2}$ \\ \toprule
	Parameter & $\dot{\bar{\lambda}}_{\text{p}, \, \text{max}}$ & $\dot{\bar{\lambda}}_{\text{p}, \, \text{min}}$ & $\mu_\text{a}$ & $\kappa$ & $\beta_1$ & $k_{2/5, \, \text{start}}$\\ \hline
	Value & $0.0000699\,\text{s}^{-1}$ & $-0.0002323\,\text{s}^{-1}$ & $11.857\,\text{kPa}$ & $148.262\,\text{kPa}$ & $0.001006\,\text{s}^{-1}$ & $1.82758\,\text{s}^{-1}$\\ \bottomrule
	\end{tabular*}
\end{table}

\Cref{fig:Arterial Wall and PharmacoMechanical Results} (right) shows the behavior of the artery with and without the influence of calcium channel blockers. In the simulation without medicinal effects, we see that the radius of the artery decreases with time and, in the end, is smaller than the radius of the initial unloaded artery. This behavior is consistent with the results of \cite{uhlmann2023chemo}. In the simulation with calcium channel-blocker influence, we observe, as expected, that the radius of the artery decreases at a much slower rate, owing to the lower $\textrm{Ca}^{2+}$ availability. Additionally, we see that with increasing drug concentration, the arterial wall softens. This can be observed in \Cref{fig:Arterial Wall and PharmacoMechanical Results} (right) in the heartbeat phase of the simulation, where the magnitude of radius change between the systolic and diastolic pressure is higher in the simulation with drug influence. The vasodilatory action of calcium channel blockers is satisfactorily captured by the model.

\begin{figure}[tb] 
	\begin{subfigure}[b]{0.5\textwidth}
		\centering
		\includegraphics[scale=0.8]{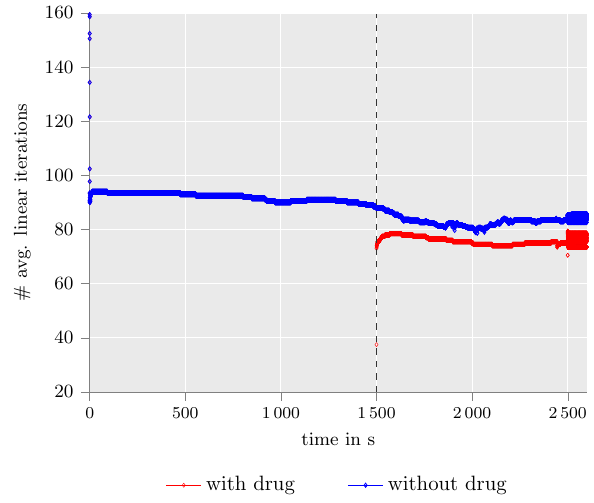}
	\end{subfigure}
		\begin{subfigure}[b]{0.5\textwidth}
		\centering
		\includegraphics[scale=0.8]{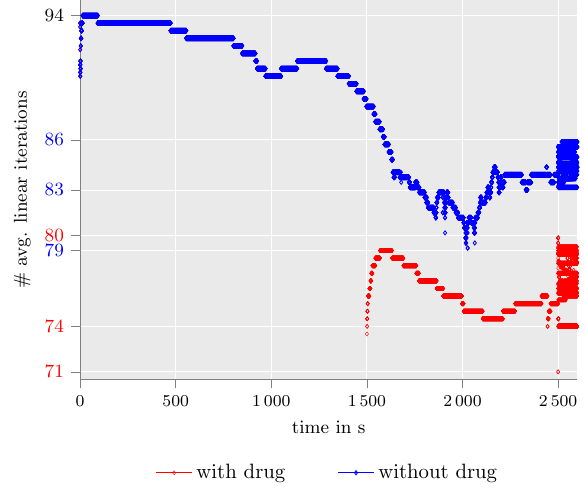}
	\end{subfigure}
	\caption{%
		Left and Right: Average number of linear iterations of loading protocol in \Cref{fig:Loading Protocol long term} and geometry in \Cref{fig:Arterial Wall and PharmacoMechanical Results} (left). 
		Drug influence is switched off for one of the simulations. 
		Computations on 50~processor cores with 50\,000 degrees of freedom. 
		Two-level overlapping Schwarz preconditioner with RGDSW coarse space and (SF) recycling. 
		Before 1\,500\,s, the iteration counts for the simulation with the influence of drugs are identical to the one without drugs. 
		Right: Detailed view between 70 and 95 average number of iterations. 
		The grid lines of the $y$ axis mark rounded maxima and minima of the iteration count with respect to the sections $[1\,500,2\,500]$ and $[2\,500,2\,600]$. 
		We observe a range of approximately 20\% (between 74 and 94) in iteration counts for simulations with and approximately 15\% (between 79 and 94) without drug influence, ignoring outliers in the loading phase. 
		The range of iteration counts increases in the heart-beat phase of the experiment, mainly due to the increased deformation caused by the softening of the material. 
		With and without the influence of drugs, respectively, we observe 74 to 80 and 83 to 86 average number of linear iterations.
	\label{fig:PharmacoMechanical Results linear iterations}}
\end{figure}

For further numerical tests and to optimize our solver, the average number of linear iterations per time step throughout the simulation are relevant and are shown in \Cref{fig:PharmacoMechanical Results linear iterations}. 
Except for higher iteration counts in the loading phase from 0\,s to 1\,s and one outlier when the diffusion process starts, the iteration count is fairly constant and even decreases over time. 
Considering this tendency, it suffices to analyze not the whole experiment but an abridged one. 
We note that, when enabling the diffusion process at $t = 1\,500$\,s, we observe a significant reduction in the number of linear iterations and, at the same time, an extreme increase in the magnitude of the nonlinear residual $\| F_k \|_{l^2}$ by more than $14$~orders of magnitude; see also~\Cref{fig:Arterial segment strong scaling linear iter}. 
We presume that there is a relation between those two effects, which we plan to further investigate in future work.

\subsection{Weak parallel scalability}
\label{sec:weak scaling}

\begin{figure}[tb]
\centering
 \begin{subfigure}[b]{0.45\textwidth}
 \centerline{\includegraphics[scale=2.7]{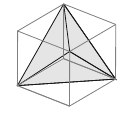}}
 \end{subfigure}
 \begin{subfigure}[b]{0.45\textwidth}
 \centerline{\includegraphics[scale=0.17]{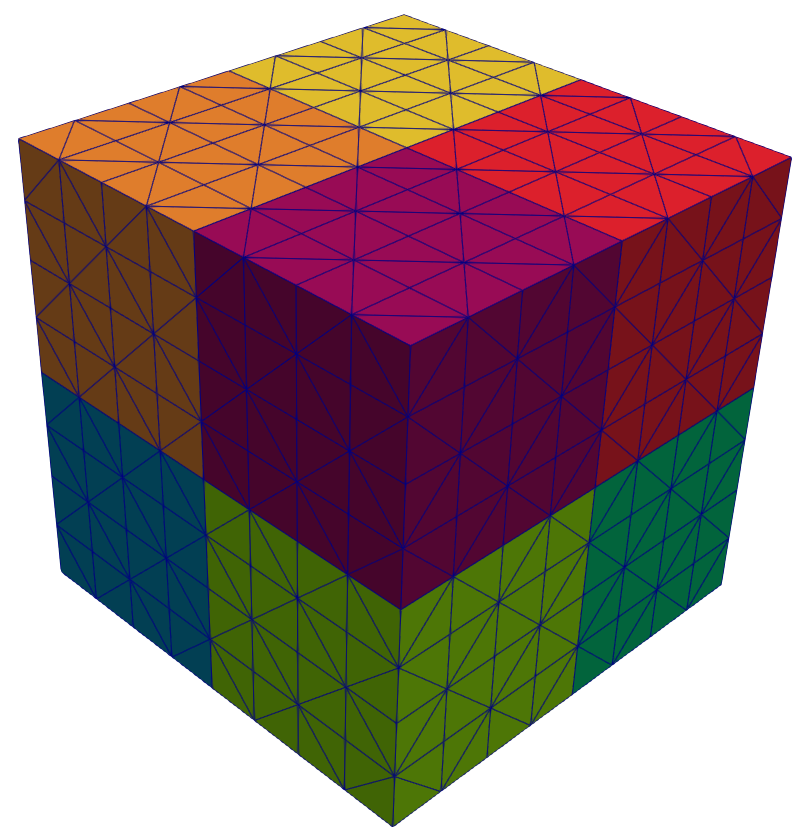}}
 \vspace{0.4cm}
 \end{subfigure}
 \caption{Partition of a cube into 5~tetrahedra (left) and corresponding structured partition of the unit cube into 8~subdomains, each containing $5{\cdot}4^3$ tetrahedra (right). 
 The orientation of a partition into 5 tetrahedra needs to be switched from one cube to a neighboring cube such that a checkerboard pattern is obtained. This pattern is visible on the right by looking at the finite element edges.}
	\label{fig:structured partition}
\end{figure}

To analyze numerical and parallel scalability, we employ structured meshes for which the number of degrees of freedom per subdomain remains constant, but the number of subdomains increases. 
For perfect scalability, the number of iterations and time-to-solution should stay asymptotically constant with an increase in the number of subdomains and processor cores. 
A voxel mesh of the unit cube is first created and subsequently refined with 5~tetrahedra per voxel; see \Cref{fig:structured partition} (left). 
The mesh of the cube is then partitioned into subcubes as depicted in \Cref{fig:structured partition} (right). 
We choose $5{\cdot}6^3$ tetrahedra per subdomain, which results in 7\,924 degrees of freedom per nonoverlapping subdomain. For the overlapping subdomains, we have a minimum of 11\,884 and maximum of 16\,7516 degrees of freedom per subdomain. 
The number of subdomains or processor cores is increased from 216 to 1\,000 to 4\,096, which produces systems with approximately 1.3, 6.2, and 25.2 million degrees of freedom, respectively.

\begin{figure}[tb] 
\begin{subfigure}[b]{0.4\textwidth}
	 	\centering
	 	\includegraphics[scale=1.0]{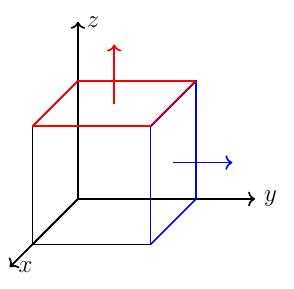}
	 \end{subfigure}
 	 \begin{subfigure}[b]{0.6\textwidth}
 	 	\centering
		\includegraphics[scale=1.15]{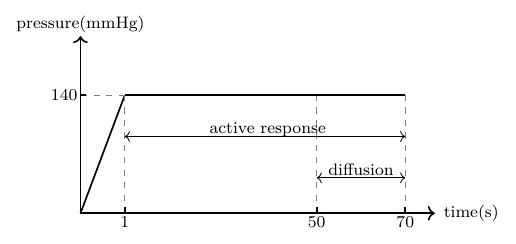}
	 \end{subfigure}
	 
	\caption{Setup for weak scaling tests on unit cube (side length: 1\,mm). 
	Left: Pressure load on $z=1$ plane and $y=1$ plane. 
	Pulling in normal directions of the respective planes. 
	Boundary conditions for the displacement: Fixed at $(0, 0, 0)$ in all directions. 
	The three faces adjacent to $(0, 0, 0)$ are fixed in their respective normal directions. 
	Right: Loading protocol for weak scaling with diffusion process and active response. Loading from 0\,s to 1\,s with 20 load steps. Active response from 1\,s until the end of 70\,s  with $\Delta t=0.25$\,s. Diffusion starts at $t=50$\,s.
	\label{fig:Cube setup weak scaling}}
\end{figure}

We consider a simplified setup of an axial pulling of the unit cube (\Cref{fig:Cube setup weak scaling}, left) with the presented material model from \Cref{sec:modeling} and model settings from \Cref{tab:passive,tab:chem_nOpt,tab:opt}. 
The unit cube has a side length of 1\,mm. 
The boundary conditions for the displacement are prescribed as follows: 
At $(0,0,0)$, the cube is fixed in all directions. 
The three faces adjacent to $(0,0,0)$ are fixed in their respective normal directions. 
The cube is then pulled at the top face, corresponding to the $z=1$ plane, in $z$~direction and at the right face, corresponding to the $y=1$ plane, in $y$~direction. 
For diffusion, a Dirichlet boundary condition emulating drug inflow is prescribed at the back face, corresponding to the $x=0$ plane.

For the weak scaling tests, we choose a relatively simple loading protocol as shown in \Cref{fig:Cube setup weak scaling}. 
From 0\,s to 1\,s, a load is built up and then kept constant at 140\,mmHg. 
Compared to \Cref{sec:Pharmaco-mechanicalResults}, we choose a larger diffusion coefficient of $D=6 \cdot 10^{-3}$ to increase the influence of the diffusion, which is diminished due to the shorter simulation time.

For this setup, we will test different coarse spaces (cf.~\Cref{subsection:GDSWandRGDSW}) and recycling strategies as defined in \Cref{sec:Recycling Strategies}. 
Generally, we can expect that with more extensive reuse of coarse space information, the setup cost of the preconditioner will decrease while the linear iteration count and, thus, the runtime of the solver will increase. 
Therefore, we need to consider carefully which configuration yields the most time gain.

Using the GDSW coarse space, the method does not scale well (cf.~\Cref{tab:GDSW vs RGDSW:weak scaling}): 
this is mainly a result of increasing setup times for an increasing number of subdomains since the GDSW coarse space is too large to be solved efficiently by a direct solver. 
Notably, with 4\,096 subdomains, the computational time exceeds a 24~hour limit at the Fritz supercomputer. 
A three-level variant may accelerate the solution; cf.~\cite{heinlein:2023:pst,heinlein:2020:tle}. 
The weak parallel scalability is significantly better for the RGDSW coarse space, since the smaller coarse space dimension (cf. \Cref{tab:WeakScaling Coarse Space Dimension}) reduces the time, for instance, for the factorization of the coarse matrix and the application of the factorization.

\begin{table}[tb]
\caption{Coarse space dimension for the two-level monolithic overlapping Schwarz preconditioners with a GDSW and RGDSW coarse space for the structure-chemical interaction system and different number of cores. Geometry in \Cref{fig:Cube setup weak scaling}.}
\label{tab:WeakScaling Coarse Space Dimension}
\begin{tabular*}{\textwidth}{@{\extracolsep\fill}lrrr@{\extracolsep\fill}}
\toprule
\textbf{Coarse Space} &  \textbf{216} & \textbf{1\,000} & \textbf{4\,096} \\
\midrule
GDSW        					   &  7\,430 &  38\,826 & 169\,740   \\
RGDSW      					   &  875 &  5\,103 & 23\,625   \\
\bottomrule
\end{tabular*}
\end{table}

Also for RGDSW, the parallel efficiency obtained for this complex nonlinear coupled problem is significantly below the efficiency routinely reported for domain decomposition applied to single physics problems. This is to some extent a result of an increase in GMRES iterations within the range of processor cores considered here. Specifically, when increasing the number of processor cores, we notice a significant increase in the orthogonalization time of GMRES: for instance, for the (SF)+(CM) combination of recycling strategies (cf. \Cref{tab:RGDSW:weak scaling and recycling}), we obtain $486$\,s, $1\,136$\,s, and $1\,978$\,s for $216$, $1\,000$, and $4\,096$ cores, respectively. The significant increase in total time cannot be exclusively explained by the increasing number of GMRES iterations. As we do not observe this effect in our strong scalability studies in~\Cref{sec:strong scaling}, where the increase in GMRES iterations is even larger, we suspect that this effect also relates to increasing Tpetra communication times. However, the main part of the increase in the solve time can be attributed to the increasing dimension of the coarse space respectively the coarse matrix. Despite the significant reduction compared to the classical GDSW coarse space, the coarse space dimension for RGDSW still increases from $875$ for $216$ cores to $23\,625$ for $4\,096$ cores. Whether numerical scalability is obtained for a larger number of subdomains and cores remains to be tested. A hotspot analysis using more detailed subtimers should be performed in the future.

In \Cref{tab:GDSW:recycling,tab:RGDSW:weak scaling and recycling}, we compare the different recycling strategies introduced in \Cref{sec:Recycling Strategies}. 
According to these results, the (SF) strategy performs best. 
It decreases the setup time significantly by 35--45\% compared to no use of recycling, and the average iteration counts do not increase. 
The reuse of the coarse basis (CB) and the coarse matrix (CM) are much more invasive strategies. 
Nevertheless, in the considered examples, they outperform the no-reuse strategy for the RGDSW coarse space. 
For the GDSW coarse space, the (CM)+(SF) strategy performed worst with respect to iteration count and solver time. 
With the RGDSW coarse space, the (CB)+(SF) strategy delivered the worst results with respect to iteration counts and solver time. 
The results for the (CM)+(SF) strategy with the RGDSW coarse space are comparable to the (SF) strategy as it saves further time in the preconditioner setup, even though the iteration count increases. 

\begin{table}[tb]
\caption{Results for different recycling strategies (cf.~\Cref{sec:Recycling Strategies}) and the two-level monolithic preconditioners with a GDSW coarse space for the structure-chemical interaction system with an overlap of $\widehat{\delta}=1$. Simulation of 296 time steps. Loading protocol and geometry in \Cref{fig:Cube setup weak scaling}. The best total time is highlighted in boldface.}
\label{tab:GDSW:recycling}
\begin{tabular*}{\textwidth}{@{\extracolsep\fill}lcrrrr@{\extracolsep\fill}}
\toprule
\multicolumn{6}{c}{\textbf{GDSW Coarse space}} \\
\midrule
\textbf{ \#\,Cores}
  & \textbf{Recycling Strategy}
	& \textbf{Avg. \# Its.}
	& \textbf{Setup$^{*}$}
	& \textbf{Solve$^{+}$}
	& \textbf{Total}
	\\
\midrule
\multirow{4}{*}{216} & --    & 123.7 (2.1) & 1\,931\,s & 2\,001\,s & 3\,932\,s \\
										 & SF    & 123.7 (2.1) & 1\,138\,s & 1\,839\,s & \textbf{2\,997\,s} \\
										 & SF+CB & 176.5 (2.1) &    811\,s & 2\,594\,s & 3\,405\,s \\
										 & SF+CM & 201.8 (2.1) &    805\,s & 3\,542\,s & 4\,347\,s \\
\midrule
\multicolumn{6}{c}{\makebox[0pt]{\scriptsize{$^{\dagger}$\textbf{avg. \# its.} of linear solver (nonlinear solver), $^{*}$\textbf{setup} of preconditioner, $^{+}$linear \textbf{solver} time}}} \\
\bottomrule
\end{tabular*}
\end{table}

\begin{table}[tb]
\caption{Weak scaling results for two-level monolithic preconditioners with GDSW and RGDSW coarse space for structure-chemical interaction system with overlap $\widehat{\delta}=1$. Simulation of 296 time steps. Loading protocol and geometry in \Cref{fig:Cube setup weak scaling}. The best total time is highlighted in boldface.}
\label{tab:GDSW vs RGDSW:weak scaling}
\begin{tabular*}{\textwidth}{@{\extracolsep\fill}l l cc rrr@{\extracolsep\fill}}
\toprule
\textbf{ \#\,Cores}
  & \textbf{Coarse Space}
  & \textbf{Recycling Strategy}
	& \textbf{Avg. \# Its.$^{\dagger}$}
	& \textbf{Setup$^{*}$}
	& \textbf{Solve$^{+}$}
	& \textbf{Total}
	\\
\midrule
\multirow{2}{*}{216}  &  GDSW & SF & 123.7 (2.1) & 1\,138\,s & 1\,839\,s & 2\,997\,s \\
                      & RGDSW & SF & 152.5 (2.1) &    789\,s & 1\,827\,s & \textbf{2\,616\,s} \\
\midrule
\multirow{2}{*}{1\,000} &  GDSW & SF & 164.1 (2.1) & 2\,414\,s & 6\,297\,s & 8\,711\,s \\
                      & RGDSW & SF & 184.7 (2.1) & 1\,006\,s & 3\,449\,s & \textbf{4\,435\,s} \\
\midrule
\multicolumn{7}{c}{\makebox[0pt]{\scriptsize{$^{\dagger}$\textbf{avg. \# its.} of linear solver (nonlinear solver), $^{*}$\textbf{setup} of preconditioner, $^{+}$linear \textbf{solver} time}}} \\
\bottomrule
\end{tabular*}
\end{table}

\begin{table}[tb]
\caption{Weak scaling results for two-level monolithic preconditioners with RGDSW coarse space for structure-chemical interaction system with overlap $\widehat{\delta}=1$. Simulation of 296 time steps. 
Loading protocol and geometry in \Cref{fig:Cube setup weak scaling}.
The best total time is highlighted in boldface.}
\label{tab:RGDSW:weak scaling and recycling}
\begin{tabular*}{\textwidth}{@{\extracolsep\fill}l crrrr@{\extracolsep\fill}}
\toprule				
\multicolumn{6}{c}{\textbf{RGDSW Coarse space}} \\
\midrule
\textbf{ \#\,Cores}
  & \textbf{Recycling Strategy}
	& \textbf{Avg. \# Its.$^{\dagger}$}
	& \textbf{Setup$^{*}$}
	& \textbf{Solve$^{+}$}
	& \textbf{Total}
	\\
\midrule

\multirow{4}{*}{216}  & --		& 152.5 (2.1)	& 1\,207\,s	& 1\,990\,s	&  3\,197\,s \\
										  & SF		& 152.5 (2.1)	&    788\,s	& 1\,827\,s	& \textbf{2\,616\,s} \\
										  & SF+CB	& 193.7 (2.1)	&    562\,s	& 2\,362\,s	&  2\,924\,s \\
										  & SF+CM & 183.8 (2.1)	&    706\,s	& 2\,096\,s	&  2\,802\,s \\
\hline
\multirow{4}{*}{1\,000} & --		& 184.7 (2.1)	& 1\,841\,s	& 4\,331\,s	&  6\,172\,s \\
										  & SF		& 184.7 (2.1)	& 1\,006\,s	& 3\,449\,s	& \textbf{4\,435\,s} \\
										  & SF+CB	& 239.0 (2.1)	&    794\,s	& 5\,049\,s	&  5\,843\,s \\
										  & SF+CM & 210.2 (2.1)	&    824\,s	& 3\,976\,s	&  4\,799\,s \\
\hline
\multirow{4}{*}{4\,096} & --	  & 219.9 (2.1)	& 2\,736\,s	& 8\,844\,s	& 11\,580\,s \\
										  & SF		& 219.9 (2.1)	& 1\,483\,s	& 6\,934\,s	& \textbf{8\,417\,s} \\
										  & SF+CB	& 286.5 (2.1)	& 1\,217\,s	& 8\,827\,s	& 10\,044\,s \\
										  & SF+CM & 245.6 (2.1)	&    860\,s	& 7\,807\,s	&  8\,666\,s \\
\midrule
\multicolumn{6}{c}{\makebox[0pt]{\scriptsize{$^{\dagger}$\textbf{avg. \# its.} of linear solver (nonlinear solver), $^{*}$\textbf{setup} of preconditioner, $^{+}$linear \textbf{solver} time}}} \\
\bottomrule
\end{tabular*}
\end{table}

\subsection{Strong parallel scalability} \label{sec:strong scaling}

\begin{figure}[tb]
	 \begin{subfigure}[b]{0.5\textwidth}
	 	\centering
	 	\includegraphics[scale=0.21]{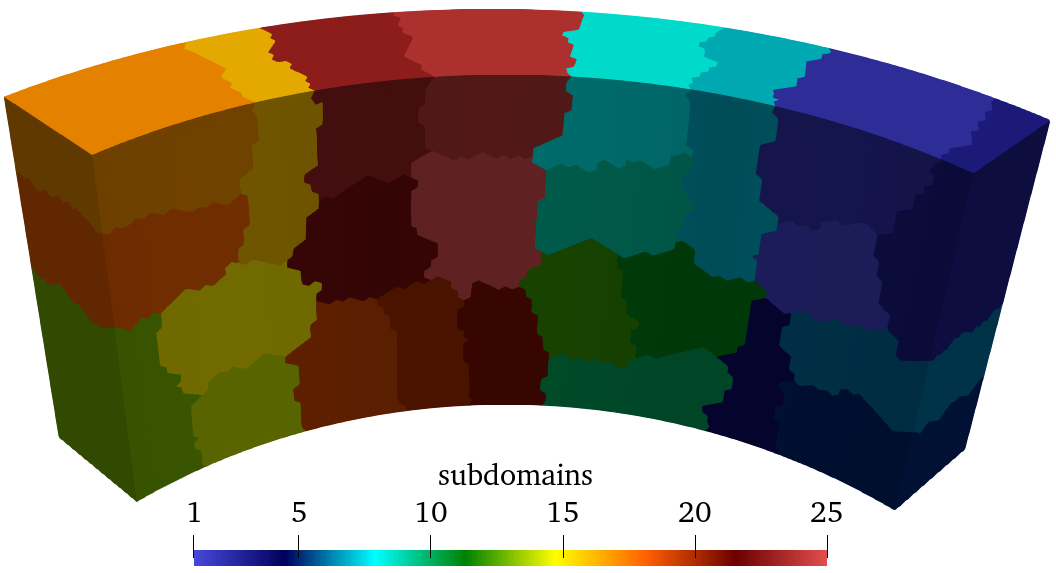}
	 		 \end{subfigure}
	 		 \begin{subfigure}[b]{0.5\textwidth}
	 		 	 	\centering
	   \includegraphics[scale=0.21]{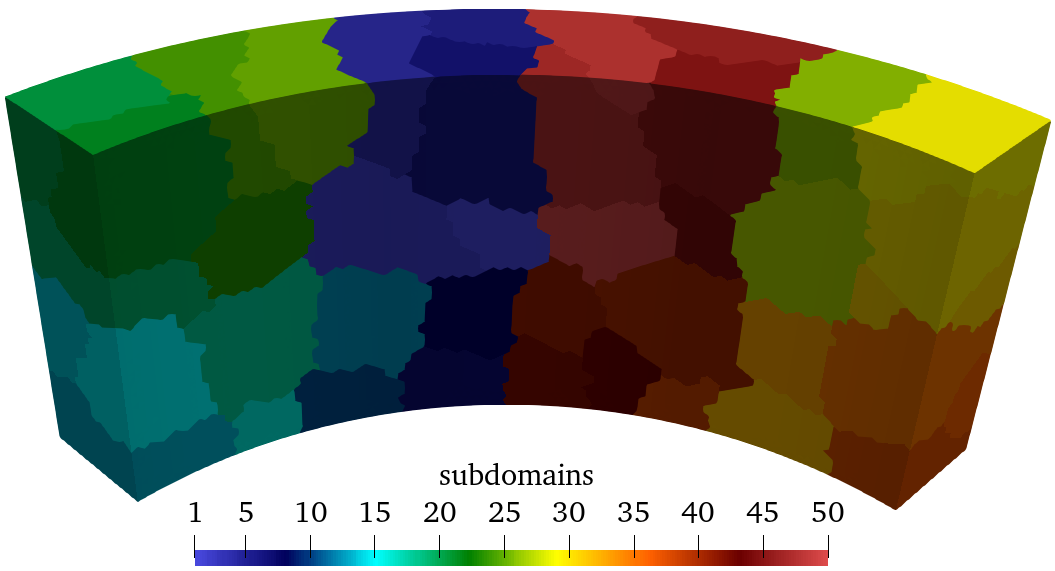}
	 \end{subfigure}
	\caption{%
		Unstructured partition of mesh with METIS into 25 and 50 subdomains. 
		Geometry in \Cref{fig:Arterial Wall and PharmacoMechanical Results}. 
		With 25~subdomains (left), the partition is almost two-dimensional with respect to the interior arterial wall. 
		With 50~subdomains (right), we observe a slight partition of the domain in the radial direction.
	\label{fig:Arterial segment partition METIS}}
\end{figure}

\begin{figure}[tb] 
\centering
		\includegraphics[scale=0.9]{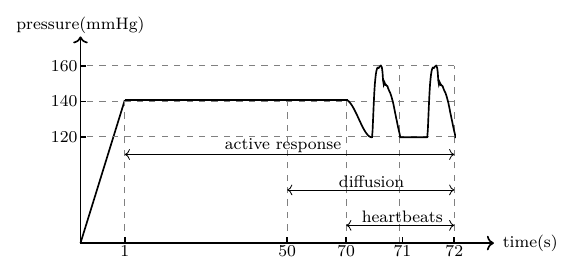}
	\caption{Loading protocol for strong scaling tests with active response, diffusion process, and two heartbeats. Loading from $t=0$\,s to $t=1$\,s with 40 load steps. A constant pressure of 140\,mmHg is applied until $t=70$\,s with $\Delta t =0.25$\,s. Starting at $t=70$\,s two heartbeats are simulated with $\Delta t=0.02$\,s. The time-dependent change of the loading during the heartbeats was generated based on the flow rate profile from~\cite{balzani:2015:nmf}, which is based on data from~\cite{hemolab}. \label{fig:Arterial segment setup strong scaling}}
\end{figure}

\begin{figure}[tb]
	\begin{subfigure}[b]{0.495\textwidth}
		\centering
		\includegraphics[scale=0.74]{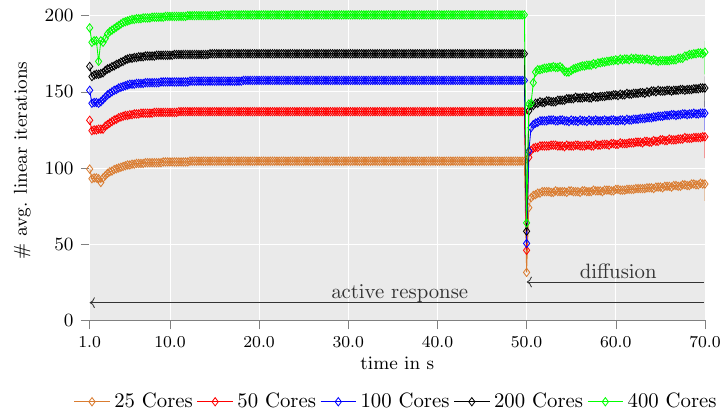}
	\end{subfigure}
	\begin{subfigure}[b]{0.505\textwidth}
		\centering
		\includegraphics[scale=0.74]{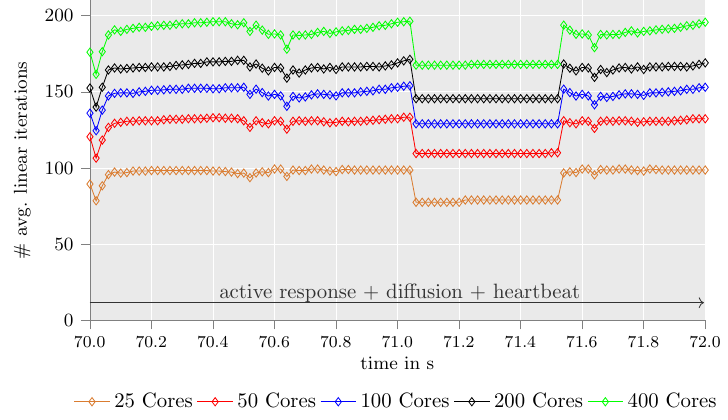}
	\end{subfigure}
	\caption{Iteration counts for different phases of loading as depicted in \Cref{fig:Arterial segment setup strong scaling}. Geometry in \Cref{fig:Arterial Wall and PharmacoMechanical Results}.  
 RGDSW coarse space with (SF) recycling strategy. 
	Compare with the results in \Cref{fig:PharmacoMechanical Results linear iterations}, where a coarser mesh was used. 
	}
	\label{fig:Arterial segment strong scaling linear iter}
\end{figure}

For the analysis of strong scaling, we keep the problem size constant while increasing the number of processor cores from~25 to~400. 
We consider the arterial segment used in \Cref{sec:Pharmaco-mechanicalResults} (cf.~\Cref{fig:Arterial Wall and PharmacoMechanical Results,fig:Arterial segment partition METIS}) and construct a mesh such that the coupled system~\cref{eq:implicit discretized system} has 2.5\,million degrees of freedom. 
The geometry is partitioned via METIS into unstructured subdomains; see \Cref{fig:Arterial segment partition METIS}. 
We apply an abridged loading protocol, similar to the weak scaling analysis.
Here, we additionally simulate two heartbeats; see \Cref{fig:Arterial segment setup strong scaling}. 

For the preconditioner, the RGDSW coarse space is employed, and we reuse the symbolic factorizations (SF) as described in \Cref{sec:Recycling Strategies}, since this configuration provided the best results for weak scaling. 
Let us note that, for a relatively small number of processor cores, GDSW provides competitive results compared to RGDSW; see the results in \Cref{tab:GDSW vs RGDSW:weak scaling} for 216~cores.

As before in \Cref{sec:Pharmaco-mechanicalResults} (cf. \Cref{fig:PharmacoMechanical Results linear iterations}), we observe a decline in the average number of linear iterations when the diffusion process starts; see also the discussion at the end of~\Cref{sec:Pharmaco-mechanicalResults}. 

\begin{figure}[tb]
	\begin{subfigure}[T]{0.495\textwidth}
		\centering
		\includegraphics[scale=0.8]{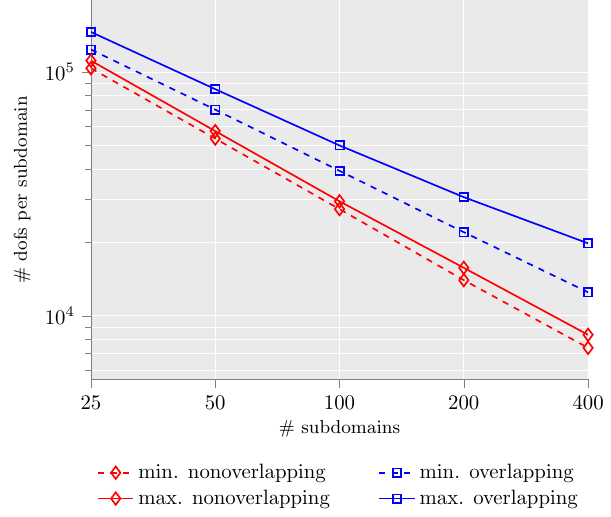}
	\end{subfigure}
	\begin{subfigure}[T]{0.495\textwidth}
		\centering	
		\includegraphics[scale=0.8]{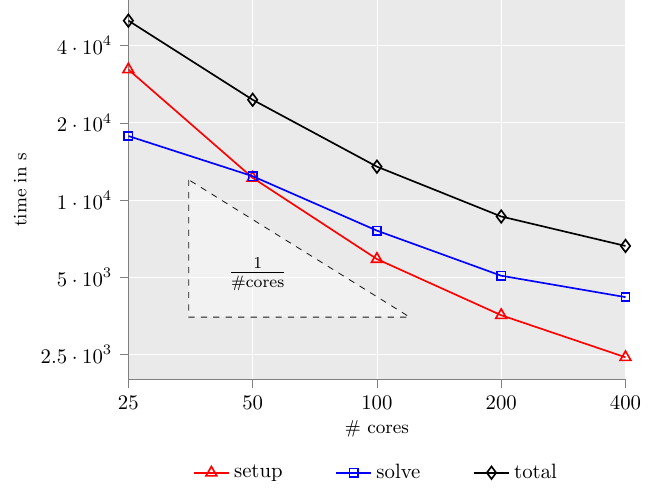}
	\end{subfigure}
	\caption{Degrees of freedom per subdomain versus number of subdomains (left). Minimum and maximum number of degrees of freedom per nonoverlapping and overlapping subdomain, respectively. 
	Unstructured mesh partition with METIS; cf. \Cref{fig:Arterial segment partition METIS}. Strong scaling plot for RGDSW with time in seconds versus number of cores for setup, solve, and total time (right); cf. \Cref{tab:RGDSW and GDSW:strong scaling}. 
	The black dashed triangle shows the slope for ideal scalability.} \label{fig:Subdomains vs. dofs}
\end{figure}

Results for GDSW and RGDSW are given in~\Cref{tab:RGDSW and GDSW:strong scaling}; the results for RGDSW are also plotted in~\Cref{fig:Subdomains vs. dofs} (right). 
For RGDSW, we initially see an efficiency of over 90\%, while it decreases considerably for a larger number of processors. 
Overall, the performance of the RDGSW coarse space is relatively robust. 
Especially the setup of the preconditioner scales well with an increasing number of subdomains. 

Before we discuss several factors that influence the performance, we consider the timings for GDSW. 
As the coarse space grows more rapidly with the number of subdomains than the RGDSW coarse space, we expect a lower efficiency. 
Indeed, the efficiency deteriorates much more quickly. 
Despite fewer iterations than using the RGDSW coarse space for 400~subdomains, the time taken by the linear solver is 10\,834\,s larger; 
we conclude that the coarse solve takes up a large portion of the linear solver time. 
In total, 13\,590\,s are spent in 186\,200 applications of the preconditioner.

\begin{table}[tb]
\caption{Strong scaling results for two-level monolithic overlapping Schwarz preconditioners with RGDSW coarse space and overlap $\widehat{\delta}=1$.
Average linear iteration count over all time steps with nonlinear solver iterations in parenthesis. Unstructured mesh partition with METIS. Simulation of 416~time steps. Loading protocol in \Cref{fig:Arterial segment setup strong scaling}, geometry in \Cref{fig:Arterial Wall and PharmacoMechanical Results}.}
\label{tab:RGDSW and GDSW:strong scaling}
\centering
\begin{tabular*}{\textwidth}{@{\extracolsep\fill}lrrrr@{\extracolsep\fill}}
\toprule
\textbf{ \#\,Cores}
	& \textbf{Avg. \# Its.$^{\dagger}$}
	& \textbf{Setup$^{*}$}
	& \textbf{Solve$^{+}$}
	& \textbf{Total}
	\\
\midrule

\multicolumn{5}{c}{\textbf{RGDSW Coarse Space}} \\
\midrule

\phantom{0}25 &  97.9 (2.4) & 32\,270\,s & 17\,770\,s & 50\,040\,s \\
\phantom{0}50 & 131.8 (2.4) & 12\,210\,s & 12\,390\,s & 24\,600\,s \\
          100 & 151.6 (2.4) &  5\,890\,s &  7\,603\,s & 13\,498\,s \\
          200 & 169.1 (2.4) &  3\,559\,s &  5\,072\,s &  8\,631\,s \\
          400 & 193.3 (2.4) &  2\,440\,s &  4\,186\,s &  6\,626\,s \\
\midrule

\multicolumn{5}{c}{\textbf{GDSW Coarse Space}} \\ 	
\midrule

\phantom{0}25 &   79.6 (2.4) & 35\,570\,s & 14\,850\,s & 50\,420\,s \\
\phantom{0}50 &  119.8 (2.4) & 20\,550\,s & 12\,590\,s & 33\,140\,s \\
          100 &  143.5 (2.4) & 12\,050\,s &  8\,918\,s & 20\,968\,s \\
          200 &  159.9 (2.4) &  9\,650\,s &  7\,406\,s & 17\,056\,s \\
          400 &  185.5 (2.4) & 11\,370\,s & 15\,020\,s & 26\,390\,s \\
\midrule
\multicolumn{5}{c}{\makebox[0pt]{\scriptsize{$^{\dagger}$\textbf{avg. \# its.} of linear solver (nonlinear solver), $^{*}$\textbf{setup} of preconditioner, $^{+}$linear \textbf{solver} time}}} \\
\bottomrule
\end{tabular*}
\end{table}

There are several properties that influence the performance of the two-level additive Schwarz method. 
For example, if the size of local problems decreases, generally, the number of iterations increases, which is also evident from \Cref{tab:RGDSW and GDSW:strong scaling}. 
A larger number of iterations not only requires more applications of the system matrix in GMRES but also of the preconditioner. 
Moreover, the complexity of the orthogonalization routine within GMRES increases quadratically with the number of iterations. 
Similarly, since the coarse space size increases with the number of subdomains (cf. \Cref{tab:WeakScaling Coarse Space Dimension}), the (at worst) cubic complexity of a direct solver results in significantly longer runtimes, which can, to some degree, be mitigated by using a parallel sparse direct solver. 
On the other hand, the subdomain problems decrease in size and, thus, profit from the superlinear complexity of direct solvers. 
Furthermore, there are hardware-related factors like cache effects, which we will not elaborate here.

Instead, we discuss how the mesh partition can influence the performance. 
The partition of the mesh into 25~subdomains is almost two-dimensional in the sense that there is usually only a single subdomain in the radial direction; see \Cref{fig:Arterial segment partition METIS} (left). 
This effect is diminished with a partition into 50~subdomains (see \Cref{fig:Arterial segment partition METIS} (right)); 
it significantly increases coupling between subdomains for the largest test case of 400~subdomains. 
The additional coupling has an influence on the performance and weakens strong scalability. 

For good load balancing, subdomains should ideally have the same number of nodes. 
Otherwise, results for strong scaling may be impacted by idle cores. 
As we can see in \Cref{fig:Subdomains vs. dofs} (left), the METIS-computed mesh partition for the nonoverlapping subdomains differs only slightly with respect to the minimum (red dashed line) and maximum (red solid line) subdomain size. 
However, locally, we do not solve nonoverlapping but overlapping subdomain problems. 
\Cref{fig:Subdomains vs. dofs} (left) shows that the minimum (blue dashed line) and maximum (blue solid line) overlapping subdomain sizes diverge with an increasing number of subdomains. 
A reason may be the shape of the subdomains. 
Even if the nonoverlapping subdomains are identical in size, a different shape can significantly influence the size of the corresponding overlapping subdomain.
For example, the ratio of nodes in the overlapping to the nodes in the nonoverlapping subdomain is large for a thin beam and small for a cube or sphere. 

This is not the only source that can lead to a deterioration in performance. 
With an increasing number of subdomains (for a fixed mesh), the total number of nodes in overlapping subdomains increases as well. 
As a result, we cannot obtain perfect scalability. 
We consider, for example, a cube as a nonoverlapping subdomain with $10{\times} 10{\times} 10$ nodes. 
If we prescribe an overlap of one layer of additional nodes, we obtain $12{\times} 12{\times} 12=1\,728$ nodes. 
By decomposing the cube into $2{\times} 2{\times} 2$ subcubes, we obtain overlapping subdomains of size $7{\times} 7{\times} 7=343$, in total $8{\times} 343=2\,744$ nodes, an increase of 58.8\% with respect to the initial overlapping cube. 
Consequently---disregarding other factors---we cannot obtain perfect strong scalability. 

A further reduction of the computing time could be obtained by a temporal homogenization of the diffusion process, which can also be combined with a parallel-in-time integration; cf.~\cite{frei:2022:tpt}. We might consider this in future investigations.

\section{Conclusions}\label{sec5}
Motivated by a need to better understand pharmaco-mechanical interactions in arteries, we have developed a computational framework to describe the effect of calcium channel blockers on the mechanical response. We adopted the material model of the arterial wall from Uhlmann and Balzani \cite{uhlmann2023chemo} and extended it to include pharmacological effects. The transmural drug transport was modeled using the reaction-diffusion equation and the resulting coupled system of equations was discretized using the finite element method. The resulting nonlinear finite element system is solved with a Newton method combined with a load stepping strategy where the tangent problems are solved using a parallel monolithic domain decomposition approach. Numerical results analyzing strong and weak parallel scalability for this strategy are presented showing the simulation capability of this approach. The pharmaco-mechanical model was implemented in AceGen and then transferred to FEDDLib using a newly developed interface. The numerical results show that this approach works well. Simulation results for an arterial section under physiological loading conditions show that the proposed model is able to qualitatively capture the vasodilatory effects of calcium channel blockers. An accurate simulation of patient-specific arteries and stress distributions would require additional considerations such as residual stresses, advection-based drug transport and realistic fiber distributions. Residual stresses can be included in our model using methods such as the open angle method or anisotropic growth models, whereas realistic fiber orientations may be obtained using arterial remodeling models. However, since these factors do not play a pivotal role in the pharmaco-mechanical interaction, their effects were not considered in this work and may be included in future studies.

\section*{Acknowledgments}

Financial funding from the Deutsche Forschungsgemeinschaft (DFG) through the Priority Program 2311 ``Robust coupling of continuum-biomechanical in silico models to establish active biological system models for later use in clinical applications - Co-design of modeling, numerics and usability'', project ID 465228106, is greatly appreciated.

The authors gratefully acknowledge the scientific support and HPC resources provided by the Erlangen National High Performance Computing Center (NHR@FAU) of the Friedrich-Alexander-Universität Erlangen-Nürnberg (FAU) under the NHR project k105be. NHR funding is provided by federal and Bavarian state authorities. NHR@FAU hardware is partially funded by the German Research Foundation (DFG) – 440719683.

\subsection*{Author contributions}
The order of the list of authors of this paper is alphabetical. Their contributions are as follows:\\
{\bf Daniel Balzani}: conceptualization, supervision, methodology, writing, review; {\bf Alexander Heinlein}: methodology, writing, programming, review; {\bf Axel Klawonn}: conceptualization, supervision, methodology, writing, review; {\bf Jascha Knepper}: methodology, writing, review; {\bf Sharan Nurani Ramesh}: methodology, writing, programming; {\bf Oliver Rheinbach}: conceptualization, supervision, methodology, writing, review; {\bf Lea Sa\ss manshausen}: methodology, writing, programming; {\bf Klemens Uhlmann}: writing, model formulation.

\subsection*{Conflict of interest}

The authors declare no potential conflict of interests.

\bibliography{bib}

\end{document}